\documentclass[11pt,a4paper]{article}
\usepackage{amssymb,amsmath}
\usepackage{color}
\usepackage{a4wide}
\newcommand {\eq} [1] {\begin{equation}\label{#1}}
\newcommand {\en} {\end{equation}}
\newcommand {\eqn}  {\begin{eqnarray}}
\newcommand {\enn} {\end{eqnarray}} 
\newcommand {\bstar}  {\begin{eqnarray*}}
\newcommand {\estar} {\end{eqnarray*}} 

\newcommand {\B} {{\bf B}}
\newcommand {\C} {{\bf C}}
\newcommand {\D} {{\bf D}}
\newcommand {\E} {{\bf E}}

\newcommand {\G} {{\bf G}}
\newcommand {\I} {{\bf I}}

\newcommand {\Q} {{\bf Q}}

\newcommand {\V} {{\bf V}}

\newcommand {\X} {{\bf X}}

\newcommand {\Z} {{\bf Z}}

\newcommand {\proof} {\par{\it Proof}. \ignorespaces}
\newcommand {\eproof}
      {\space
 {\ \vbox{\hrule\hbox{\vrule height1.3ex\hskip0.8ex\vrule}\hrule}}
 \par}
\newcommand {\setC}  {{\mathbb C}}

\newcommand {\mat} [1] {\left[\begin{array}{#1}}
\newcommand {\rix}     {\end{array}\right]}
\newtheorem{theorem}           {Theorem}

\newtheorem{corollary}[theorem]{Corollary}
\newtheorem{example}           {Example}
\newtheorem{remark}            {Remark}

\newcommand {\diag}  {\mathop{\rm diag}\nolimits}
\newcommand {\rank}  {\mathop{\rm rank}\nolimits}

\newcommand {\range} {\mathop{\rm range}\nolimits}

\newcommand {\wh} {\widehat}
\newcommand {\wt} {\widetilde}

\begin{document}
\title{Invariant subspace perturbation of a matrix with Jordan blocks}
\author{Hongguo Xu\thanks{Department of Mathematics, University of Kansas, Lawrence, KS 66045, USA.
email: {\tt feng@ku.edu}. }
}
\date{}
\maketitle
\begin{abstract} 
We investigate how invariant subspaces corresponding to a single eigenvalue will change when a matrix is perturbed. We focus on the invariant subspaces corresponding to an eigenvalue associated with the Jordan blocks that have the same size. 
An invariant subspace can be expressed as the range
of a full column matrix.  We characterize the perturbed
invariant subspaces with such matrices expressed in a sum form 
that exhibits the fractional orders. We also provide the formulas for the 
coefficient matrices associated with the zero and first fractional orders. 
The results extend  the standard invariant subspace perturbation theory.
\end{abstract}
\noindent
{\bf Keywords} Invariant subspace, Perturbation, Jordan block, Eigenvalue, Eigenvector

\noindent
{\bf AMS subject classification} 15A18,  47A55, 65F15
\section{Introduction} Eigenstructure perturbation is one of the fundamental problems in
numerical linear algebra and matrix theory. General perturbation theory is well established, e.g.,
\cite{GolV13,SteS90}.  
For eigenvalues and eigenvectors, most perturbation results are established for 
simple or semi-simple eigenvalues. Invariant and deflating subspace perturbation results are developed normally based on the assumption that  the associated eigenvalues are different from the rest. With the presence of Jordan blocks, the perturbation problem becomes  complicated. In \cite{Lid65,VisL60},  eigenvalue and eigenvector perturbation results are derived in terms of fractional orders.
Roughly speaking, if a Jordan block of size $\rho$ is slightly perturbed with a matrix of magnitude $t$, the eigenvalues of the perturbed matrix
are the original eigenvalue plus an order $t^{1/\rho}$ perturbation.
Other results can be found in \cite{Bau85,Kat80,MorBO97}.

Our goal is to establish the perturbation theory for the invariant subspaces of a matrix 
corresponding to a single eigenvalue that may not be semi-simple. 
The results will be derived by using the same 
change of variable and scaling 
techniques shown in \cite{Lid65,MorBO97}. We only consider the generic case 
that will be explained later. The results generalize the existing eigenvector perturbation theory given in \cite{Lid65}. 
Moreover, the results show how the eigenvectors and generalized eigenvectors are involved  in forming the matrices that define perturbed invariant subspaces.
Such information may be essential in other perturbation analysis. For instance, 
in \cite{MehX24} it is used for analyzing structure perturbation of  certain Hamiltonian matrices.
The whole picture of  invariant subspace perturbation is essentially obtained by combining 
the existing results \cite{SteS90} and the results presented here.

This paper is organized as follows. In the next section, for completion we provide
the eigenvalue and eigenvector perturbation results given in \cite{Lid65,MorBO97}
and introduce the basic tools used there. In Section~\ref{secmain} we present the 
main results and in Section~\ref{secfst} we provide the first fractional order invariant
subspace perturbation results. Section~\ref{con} contains our conclusions.

Throughout the paper, $\setC$ represents the complex  number field. 
$\setC^{m\times n} $ is the space of  $m\times n$ complex 
matrices, and $\setC^n$ is the $n$-dimensional complex column vector space.
For a matrix $B$, $B^\ast$  and $B^T$ are the adjoint (conjugate transpose) and transpose of $B$, respectively.
For a square matrix $A$, $\Lambda(A)$ is the spectrum of $A$. We
denote by $I_k$ and $0_k$ the $k\times k$ identity and zero matrices, respectively,
and by $0_{k\times m}$ the $k\times m$ zero matrix. For a positive number $c$, the notation $O(t^c)$ means
a term (or matrix) that converges to $0$ at the same rate as $t^c$ and
$o(t^c)$ means a term (or matrix) that converges  faster than $t^c$, both when $t\to 0$. 
We also use notations as $O\left(\mat{ccc}t^{c_{11}}&\ldots&t^{c_{1n}}\\
\vdots&\vdots&\vdots\\t^{c_{m1}}&\ldots&t^{c_{mn}}\rix\right)$ to show the
 orders of different blocks (entries) in a matrix.
For a given $r$-dimensional subspace $V\subset \setC^n$,  we call $X\in\setC^{n\times r}$ a matrix of $V$  if $\range X=V$.
\section{Preliminaries}

Let $A,D\in\setC^{n\times n}$ and $t\ge 0$ be a real parameter. Suppose $\lambda_0$
is one eigenvalue of $A$. The question in \cite{Lid65} is how the eigenvalues and their eigenvectors corresponding to $\lambda_0$ change, when $A$ is perturbed to $A+tD$.

Let $\mat{cc}\Xi&\Xi_c\rix\in\setC^{n\times n}$
be invertible and $\mat{cc}\Sigma&\Sigma_c\rix=\mat{cc}\Xi&\Xi_c\rix^{-\ast}$ such that
\eq{ta}
A \mat{cc}\Xi&\Xi_c\rix =\mat{cc}\Xi&\Xi_c\rix\mat{cc}A_{11}&0\\0&A_{22}\rix,
\en
where $A_{11}$ is the modified Jordan canonical form corresponding to the
single eigenvalue $\lambda_0$:
\eq{nila}
A_{11}=\Sigma^\ast A\Xi =\lambda_0I+N,\quad N=\diag(N_1,N_{2},\ldots,N_k),
\en
and
\[
N_j=\mat{cccc}0&I_{s_j}&&\\&0&\ddots&\\&&\ddots&I_{s_j}\\&&&0\rix_{j\times j},
\quad j=1,\ldots,k.
\]
The integers $s_1,\ldots,s_{k-1}$ are nonnegative and $s_k>0$. By convention, $N_1=0_{s_1}$. So $A$ has $s_j$ Jordan blocks of size $j$ corresponding 
to the eigenvalue $\lambda_0$, for $j=1,\ldots,k$.  If $A$ does not have Jordan blocks of size $j$ for a certain $j$,  then  $N_j$ is void ($s_j=0$). We assume
$\lambda_0\not\in \Lambda(A_{22})$.

Partition
\eq{td}
\mat{cc}\Xi&\Xi_c\rix^{-1}D\mat{cc}\Xi&\Xi_c\rix =\mat{cc}D_{11}&D_{12}\\D_{21}&D_{22}\rix,
\quad D_{11}=\Sigma^\ast D\Xi.
\en

Following the standard perturbation theory (\cite{SteS90}), when $t$ is sufficiently small,  
there is a matrix $P(t)$ that satisfies
$$
(A_{22}+tD_{22})P(t)-P(t)(A_{11}+tD_{11})
-tP(t)D_{12}P(t)+tD_{21}=0,
$$
and
\eqn
\nonumber
&&(A+tD)\mat{cc}\Xi&\Xi_c\rix\mat{cc}I&0\\P(t)&I\rix\\
\label{simad}
&&=\mat{cc}\Xi&\Xi_c\rix\mat{cc}I&0\\P(t)&I\rix\mat{cc}A_{11}
+t(D_{11}+D_{12}P(t))&tD_{12}\\0&
A_{22}+t(D_{22}-P(t)D_{12})\rix.
\enn
The matrix $P(t)$ has the series form (\cite[Section 2.1.4]{Kat80})
$
P(t)=\sum_{i=1}^\infty t^iP_i,
$
and $P_1$ is the solution of the Sylvester equation
$
A_{22}P_1-P_1A_{11}+D_{21}=0.
$
With (\ref{simad}) the original problem is reduced to the one about the block $A_{11}+tD_{11}
+D_{12}P(t))$.

We first consider the problem with the matrix 
\eq{eqmatrix}
A_{11}+tD_{11}=\lambda_0 I+N+tD_{11}.
\en
The results derived below can be easily extended to the case when $D_{11}$ is analytic (as $D_{11}+D_{12}P(t)$ given above),
therefore can be extended to the general case with $A+tD$. 

Partition 
\eqn
\nonumber
\Xi=\mat{ccc}\Xi_1&\ldots&\Xi_k\rix,\quad
\Sigma=\Xi^{-\ast}=\mat{ccc}\Sigma_1&\ldots&\Sigma_k\rix,\\
\label{bkxi}\\
\nonumber
\Xi_j=\mat{ccc}\Xi_{j1}&\ldots&\Xi_{jj}\rix,\quad \Sigma_j
=\mat{ccc}\Sigma_{j1}&\ldots&\Sigma_{jj}\rix,
\qquad j=1,\ldots,k
\enn
with the column numbers following the block sizes of $N$ and $N_1,\ldots,N_k$.
Let 
\eq{eigspace}
\wt \Xi_i:=\mat{cccc}\Xi_{i1}&\Xi_{i+1,1}&\ldots&\Xi_{k1}\rix,\quad
\wt\Sigma_i:=\mat{cccc}\Sigma_{ii}&\Sigma_{i+1,i+1}\ldots&\Sigma_{kk}\rix,
\en
for $i=1,\ldots,k$. Clearly, $\wt \Xi_1$ and $\wt \Sigma_1$
are matrices of the right and left eigenvector spaces of $A$, respectively. Denote
\eq{blkb}
D_{11}=\Sigma^\ast D\Xi =[B_{ij}]_{k\times k},\quad
B_{ij}=\mat{cccc}B^{(ij)}_{11}&\ldots&B^{(ij)}_{1j}\\\vdots&\ddots&\vdots\\
B^{(ij)}_{i1}&\ldots&B^{(ij)}_{ij}\rix,\quad
 B^{(ij)}_{\ell m} \in\setC^{s_i\times s_j},
\en
where $B_{ij}=\Sigma_i^\ast D\Xi_j$ 
 and $B^{(ij)}_{\ell m}=\Sigma_{i\ell}^\ast D\Xi_{jm}$,
for $1\le \ell\le i,\quad1\le m\le j$, and $i,j=1,\ldots,k$.
Define
\eq{defw}
W_i:=\mat{ccc}B_{i1}^{(ii)}&\ldots&B_{i1}^{(ik)}\\\vdots&\ddots&\vdots\\B_{k1}^{(ki)}&\ldots&B_{k1}^{(kk)}\rix=\wt \Sigma_i^\ast D\wt\Xi_i,\quad i=1,\ldots,k.
\en 
Suppose for a certain $\rho<k$, $W_{\rho+1}$ is invertible. 
The matrix pencil
\eq{pen}
\gamma\mat{cc}I_{s_\rho}&0\\0&0\rix-W_\rho
\en
has $s_\rho$
finite eigenvalues $\gamma_i^{(\rho)}$, $i=1,\ldots,s_\rho$. If $\rho=k$, the pencil
(\ref{pen}) is just $\gamma I_{s_k}-W_k$.

The following eigenvalue and eigenvector perturbation results were originally given
\cite{Lid65} and rephrased in \cite{MorBO97}.
%
 
\begin{theorem}\label{thm1}
Suppose  $A$ and $D$ are transformed to (\ref{ta}) and (\ref{td}), where
$A_{11}$ and $D_{11}$ are given in (\ref{nila}) and (\ref{blkb}), respectively,
and $t>0$ is a sufficiently small number.   
\begin{enumerate}
\item[(a)] 
Suppose  $\rho\in\{1,\ldots,k\}$ and  $W_{\rho+1}$ defined in (\ref{defw}) is invertible if $\rho<k$. Let $\gamma_1^{(\rho)},\ldots,\gamma_{s_\rho}^{(\rho)}$
be the $s_\rho$ finite eigenvalues of the pencil (\ref{pen}), and for each $i\in\{1,\ldots,s_\rho\}$, 
$
\mu_{i1}^{(\rho)},\ldots,\mu_{i\rho}^{(\rho)}
$
be  the $\rho$th roots of $\gamma_i^{(\rho)}$, i.e., 
$(\mu_{ij}^{(\rho)})^\rho=\gamma_i^{(\rho)}$, for $j=1,\ldots,\rho$. Then
$A+tD$ has $\rho s_\rho$ eigenvalues 
$$\lambda_{ij}^{(\rho)}=\lambda_0+t^{1/\rho}\mu_{ij}^{(\rho)}+o(t^{1/\rho}),\quad
i=1,\ldots, s_\rho,\quad j=1,\ldots,\rho.
$$
In addition, suppose  $\gamma_i^{(\rho)}$ is a simple eigenvalue for some $i\in\{1,\ldots,s_\rho\}$. Then the corresponding perturbed eigenvalues have a series form
$$
\lambda_{ij}^{(\rho)}=
\lambda_0+t^{1/\rho}\mu_{ij}^{(\rho)}+\sum_{\ell=2}^\infty \tau_{ij,\ell}^{(\rho)}t^{\ell/\rho},
\qquad j=1,\ldots,\rho.
$$
\item[(b)] Suppose all $W_1,\ldots,W_k$ are invertible and for a certain $\rho\in\{1,\ldots,k\}$
the finite eigenvalues $\gamma_1^{(\rho)},\ldots,\gamma_{s_\rho}^{(\rho)}$ of the pencil
(\ref{pen}) are distinct.
Let $\nu_i^{(\rho)}$ be an eigenvector of the pencil (\ref{pen}) corresponding to 
$\gamma_i^{(\rho)}$, for $i=1,\ldots,s_\rho$.  For each eigenvalue $\lambda_{ij}^{(\rho)}$
given in (a), $A+tD$ has a corresponding  eigenvector in a series form
$$
x_{ij}^{(\rho)}= \wt \Xi_\rho\nu_i^{(\rho)}
+\sum_{\ell=1}^\infty v_{ij,\ell}^{(\rho)}t^{\ell/\rho}, \qquad j=1,\ldots,\rho,
$$
where  $\wt\Xi_\rho$ is given in (\ref{eigspace}).
\end{enumerate}
\end{theorem}
In \cite{Lid65} it requires $\det W_\rho\ne 0$ for Theorem~\ref{thm1}(a), 
which is removed in \cite{MorBO97} allowing some $\gamma_i^{(\rho)}$ to be zero.
The {\em generic condition} that all $W_i$ are invertible shown in Theorem~\ref{thm1}(b) is to avoid the situation that eigenvalues perturbed from the Jordan blocks of size $\rho$
are mixed with those perturbed from Jordan blocks of sizes
other than $\rho$, as shown in \cite{MorBO97}. 

\medskip
Change of variable and scaling techniques are used in \cite{Lid65}
to transform the eigenvalue problem (\ref{eqmatrix}) to a generalized eigenvalue problem
that separates the eigenvalues under consideration from others,
which is the key tool for  proving Theorem~\ref{thm1}.
Since our invariant subspace perturbation results will be derived from the same
generalized eigenvalue problem,  we describe the transformation process below.

For a given $\rho\in\{1,\ldots,k\}$, introduce the variables 
$z=t^{1/\rho}$ and $\mu=(\lambda-\lambda_0)/z$.
Then the eigenvalue problem (\ref{eqmatrix}) becomes
\bstar
&&
\lambda I-(\lambda_0I+N+tD_{11})=z\mu I-(N+z^{\rho}D_{11})\\
&&
\qquad\qquad\qquad
=\mat{cccc}z\mu I-N_1-z^{\rho}B_{11}&-z^{\rho}B_{12}&\ldots&-z^{\rho}B_{1k}\\
-z^{\rho}B_{21}&z\mu I-N_{2}-z^{\rho}B_{22}&\ldots&-z^{\rho}B_{2k}\\
\vdots&\vdots&\ddots&\vdots\\
-z^{\rho}B_{k1}&-z^\rho B_{k2}&\ldots&z\mu I-N_k-z^{\rho}B_{kk}\rix.
\estar

The eigenvalue problem is transformed further with
diagonal scaling.
Define
$$
L_i=\left\{\begin{array}{ll}
\diag(\underbrace{I_{s_i},\ldots,I_{s_i}}_{i-\rho},\underbrace{z^{-1}I_{s_i},\ldots,z^{-\rho}I_{s_i}}_{\rho}),
&k\ge i> \rho\\
\diag(\underbrace{z^{-1}I_{s_i},\ldots,z^{-i}I_{s_i}}_i),&1\le i\le\rho\end{array}\right.
$$
and
$$
R_i=\left\{\begin{array}{ll} 
\diag(\underbrace{I_{s_i},\ldots,I_{s_i}}_{i-\rho}, 
\underbrace{I_{s_i},zI_{s_i},\ldots,z^{\rho-1}I_{s_i}}_\rho),
&k\ge i> \rho\\
\diag(\underbrace{I_{s_i},zI_{s_i}\ldots,z^{i-1}I_{s_i}}_i),&1\le i\le \rho.\end{array}\right.
$$
Then
$$
L_i(z^{\rho}B_{ij})R_{j}
=\left\{\begin{array}{l}\small
\mat{ccccccc}
z^{\rho}B_{11}^{(ij)}&\ldots&z^{\rho}B_{1,j-\rho}^{(ij)}&z^{\rho}B_{1,j-\rho+1}^{(ij)}
&z^{\rho+1}B_{1,j-\rho+2}^{(ij)}&\ldots&z^{2\rho-1}B_{1j}^{(ij)}\\
\vdots&\vdots&\vdots&\vdots&\vdots&\vdots&\vdots\\
z^{\rho}B_{i-\rho,1}^{(ij)}&\ldots&z^{\rho}B_{i-\rho,j-\rho}^{(ij)}
&z^{\rho}B_{i-\rho,j-\rho+1}^{(ij)}&z^{\rho+1}B_{i-\rho,j-\rho+2}^{(ij)}&\ldots
&z^{2\rho-1}B_{i-\rho,j}^{(ij)}\\
\vdots&\vdots&\vdots&\vdots&\vdots&\vdots&\vdots\\
zB^{(ij)}_{i-1,1}&\ldots&zB^{(ij)}_{i-1,j-\rho}&zB^{(ij)}_{i-1,j-\rho+1}
&z^2B_{i-1,j-\rho+2}^{(ij)}&\ldots&z^{\rho}B_{i-1,j}^{(ij)}\\
B^{(ij)}_{i1}&\ldots&B^{(ij)}_{i,j-\rho}&B^{(ij)}_{i,j-\rho+1}
&zB^{(ij)}_{i,j-\rho+2}&\ldots&z^{\rho-1}B_{ij}^{(ij)}\rix, \\ \\
\quad\hspace{11cm} k\ge i,j> \rho\\ \\
\mat{cccc}
z^{\rho}B^{(ij)}_{11}&z^{\rho+1}B^{(ij)}_{12}&\ldots&z^{\rho+j-1}B^{(ij)}_{1j}\\
\vdots&\vdots&\ldots&\vdots\\
z^{\rho}B^{(ij)}_{i-\rho,1}&z^{\rho+1}B^{(ij)}_{i-\rho,2}&\ldots&z^{\rho+j-1}B^{(ij)}_{i-\rho,j}\\
\vdots&\vdots&\ldots&\vdots\\
zB^{(ij)}_{i-1,1}&z^2B^{(ij)}_{i-1,2}&\ldots&z^{j}B^{(ij)}_{i-1,j}\\
B^{(ij)}_{i1}&zB^{(ij)}_{i2}&\ldots&z^{j-1}B^{(ij)}_{ij}\rix, \qquad k\ge i> \rho\ge j\ge 1\\ \\
\small \mat{ccccccc}
z^{\rho-1}B_{11}^{(ij)}&\ldots&z^{\rho-1}B_{1,j-\rho}^{(ij)}&z^{\rho-1}B_{1,j-\rho+1}^{(ij)}
&z^{\rho}B_{1,j-\rho+2}^{(ij)}&\ldots&z^{2\rho-2}B_{1j}^{(ij)}\\
z^{\rho-2}B_{21}^{(ij)}&\ldots&z^{\rho-2}B_{2,j-\rho}^{(ij)}&z^{\rho-2}B_{2,j-\rho+1}^{(ij)}
&z^{\rho-1}B_{2,j-\rho+2}^{(ij)}&\ldots&z^{2\rho-3}B_{2j}^{(ij)}\\
\vdots&\vdots&\vdots&\vdots&\vdots&\vdots&\vdots\\
z^{\rho-i}B_{i1}^{(ij)}&\ldots&z^{\rho-i}B_{i,j-\rho}^{(ij)}&z^{\rho-i}B_{i,j-\rho+1}^{(ij)}
&z^{\rho-i+1}B_{i,j-\rho+2}^{(ij)}&\ldots&z^{2\rho-i-1}B_{ij}^{(ij)}\rix, \\
\quad\hspace{10cm} 1\le i\le\rho < j\le k\\ \\
\mat{cccc}z^{\rho-1}B_{11}^{(ij)}&z^{\rho}B_{12}^{(ij)}&\ldots &z^{\rho+j-2}B_{1j}^{(ij)}\\
z^{\rho-2}B_{21}^{(ij)}&z^{\rho-1}B_{22}^{(ij)}&\ldots&z^{\rho+j-3}B_{2j}^{(ij)}\\
\vdots&\vdots&\vdots&\vdots\\
z^{\rho-i}B_{i1}^{(ij)}&z^{\rho-i+1}B_{i2}^{(ij)}&\ldots&z^{\rho+j-i-1}B_{ij}^{(ij)}\rix, 
\qquad1\le i,j\le  \rho,
\end{array}
\right.
$$
$$
L_i (z\mu I-N_i)R_i=
\left\{\begin{array}{ll} \mu I-N_i,&1\le i\le \rho\\ \\
\mat{ccccccc}z\mu I_{s_i}&-I_{s_i}&&&&&\\
&z\mu I_{s_i}&\ddots&&&&\\
&&\ddots&-I_{s_i}&&&\\
&&&z\mu I_{s_i}&-I_{s_i}&&\\
&&&&\mu I_{s_i}&\ddots&\\
&&&&&\ddots&-I_{s_i}\\
&&&&&&\mu I_{s_i}\rix,&k\ge i> \rho.\end{array}\right.
$$
Let 
\eq{deflr}
L=\diag(L_1,\ldots,L_k),\qquad R=\diag(R_1,\ldots,R_k).
\en 
Then
\eq{abuv}
\mu U(z)-V(z):=L(z\mu I-(N+z^{\rho}D_{11}))R=\mu(U+zE_U)-(V+E_V(z)),
\en
where
\bstar
U&=&\diag(I_{s_1};I_{2s_2};\ldots;I_{\rho s_\rho};0_{s_{\rho+1}},I_{\rho s_{\rho+1}};\ldots;
0_{(k-\rho)s_{k}},I_{\rho s_k}),\\
E_U&=&\diag(0_{s_1};0_{2s_2};\ldots;0_{\rho s_\rho};I_{s_{\rho+1}},0_{\rho s_{\rho+1}};\ldots;
I_{(k-\rho)s_{k}},0_{\rho s_k}),
\estar
{\small
$$V=
\mat{c|c|cccc|c|ccccccc}
0_{s_1}&&&&&&&&&&&&&\\\hline
&\ddots&&&&&&&&&&&&\\\hline
&&0&I_{s_\rho}&&&&&&&&&&\\
&&&0&\ddots&&&&&&&&&\\
&&&&\ddots&I_{s_\rho}&&&&&&&&\\
B_{\rho1}^{(\rho1)}&\ldots&B_{\rho 1}^{(\rho\rho)}&&&0&\ldots&B_{\rho1}^{(\rho k)}
&\ldots&\ldots&B_{\rho,k-\rho}^{(\rho k)}&B_{\rho,k-\rho+1}^{(\rho k)}&&\\\hline
&&&&&&\ddots&&&&&&&\\\hline
&&&&&&&0&I_{s_k}&&&&&\\
&&&&&&&&0&\ddots&&&&\\
&&&&&&&&&\ddots&I_{s_k}&&&\\
&&&&&&&&&&0&I_{s_k}&&\\
&&&&&&&&&&&0&\ddots&\\
&&&&&&&&&&&&\ddots&I_{s_k}\\
B_{k1}^{(k1)}&\ldots&B_{k1}^{(k\rho)}&&&&\ldots&B_{k1}^{(kk)}&\ldots&\ldots&B_{k,k-\rho}^{(kk)}&B_{k,k-\rho+1}^{(kk)}&&0
\rix,
$$
}
and
$$
E_V(z)=L(N+z^\rho D_{11})R-V=[E_{ij}^V(z)]_{k\times k}
$$
with
$$
E_{ij}^V(z)=\left\{\begin{array}{l}L_i(z^\rho B_{ij})R_j,\qquad 1\le i\le \rho-1\\ \\
\small
\mat{ccccccc}
z^{\rho}B_{11}^{(ij)}&\ldots&z^{\rho}B_{1,j-\rho}^{(ij)}&z^{\rho}B_{1,j-\rho+1}^{(ij)}
&z^{\rho+1}B_{1,j-\rho+2}^{(ij)}&\ldots&z^{2\rho-1}B_{1j}^{(ij)}\\
\vdots&\vdots&\vdots&\vdots&\vdots&\vdots&\vdots\\
z^{\rho}B_{i-\rho,1}^{(ij)}&\ldots&z^{\rho}B_{i-\rho,j-\rho}^{(ij)}
&z^{\rho}B_{i-\rho,j-\rho+1}^{(ij)}&z^{\rho+1}B_{i-\rho,j-\rho+2}^{(ij)}&\ldots
&z^{2\rho-1}B_{i-\rho,j}^{(ij)}\\
\vdots&\vdots&\vdots&\vdots&\vdots&\vdots&\vdots\\
zB^{(ij)}_{i-1,1}&\ldots&zB^{(ij)}_{i-1,j-\rho}&zB^{(ij)}_{i-1,j-\rho+1}
&z^2B_{i-1,j-\rho+2}^{(ij)}&\ldots&z^{\rho}B_{i-1,j}^{(ij)}\\
0&\ldots&0&0
&zB^{(ij)}_{i,j-\rho+2}&\ldots&z^{\rho-1}B_{ij}^{(ij)}\rix, \\ \\
\quad \hspace{11cm} k\ge i,j\ge \rho
\end{array}
\right.
$$
$$
E_{ij}^V(z)=
\mat{cccc}
z^{\rho}B^{(ij)}_{11}&z^{\rho+1}B^{(ij)}_{12}&\ldots&z^{\rho+j-1}B^{(ij)}_{1j}\\
\vdots&\vdots&\ldots&\vdots\\
z^{\rho}B^{(ij)}_{i-\rho,1}&z^{\rho+1}B^{(ij)}_{i-\rho,2}&\ldots&z^{\rho+j-1}B^{(ij)}_{i-\rho,j}\\
\vdots&\vdots&\ldots&\vdots\\
zB^{(ij)}_{i-1,1}&z^2B^{(ij)}_{i-1,2}&\ldots&z^{j}B^{(ij)}_{i-1,j}\\
0&zB^{(ij)}_{i2}&\ldots&z^{j-1}B^{(ij)}_{ij}\rix, \qquad k\ge i\ge \rho> j\ge 1.
$$
Observe that both $U(z)$ and $V(z)$ are matrix polynomials of $z$
and $E_V(z)=O(z)$.
When $z=0$, the pencil becomes $\mu U(0)-V(0)=\mu U-V$.  

To exploit the eigenstructure of $\mu U-V$, another equivalence
transformation is applied with block row and column permutations plus
block elementary 
eliminations.
Let $\Pi_L$
be the block row permutation that moves the block rows in $V$ led by
$B_{i1}^{(i1)}$, $i=\rho+1,\ldots,k,$ in that order right underneath the block row
led by $B_{\rho1}^{(\rho1)}$. Let $\Pi_R$
be the block column permutation  that moves the block columns in $V$ containing 
$B_{\rho1}^{(\rho j)}$,
$j=\rho+2,\ldots,k$, in that order to the right of the block column containing 
$B_{\rho1}^{(\rho,\rho+1)}$.
Then 
$
\Pi_L(\mu U-V)\Pi_R=
$
{\tiny
$$
\mu\mat{c|cc|c|cccc|c|cccc|c|ccccc}
I_{s_1}&&&&&&&&&&&&&&&&&&\\\hline
&I_{s_2}&&&&&&&&&&&&&&&&&\\
&&I_{s_2}&&&&&&&&&&&&&&&&\\\hline
&&&\ddots&&&&&&&&&&&&&&&\\\hline
&&&&I_{s_\rho}&&&&&&&&&&&&&&\\
&&&&&\ddots&&&&&&&&&&&&&\\
&&&&&&I_{s_\rho}&&&&&&&&&&&&\\
&&&&&&&I_{s_\rho}&&&&&&\ldots&&&&&\\\hline
&&&&&&&&0&0&&&\wh I_{\rho+1}&\ldots&0&0&&&\wh I_{k}\\\hline
&&&&&&&&&0&&&&&&&&\\
&&&&&&&&&I_{s_{\rho+1}}&\ddots&&&&&&&\\
&&&&&&&&&&\ddots&\ddots&&&&&&\\
&&&&&&&&&&&I_{s_{\rho+1}}&0&&&&&\\\hline
&&&&&&&&&&&&&\ddots&&&&&\\\hline
&&&&&&&&&&&&&&0_{(k-\rho-1)s_k}&&&&\\
&&&&&&&&&&&&&&&0&&&\\
&&&&&&&&&&&&&&&I_{s_k}&\ddots&&\\
&&&&&&&&&&&&&&&&\ddots&\ddots&\\
&&&&&&&&&&&&&&&&&I_{s_k}&0
\rix
$$
$$
-\mat{c|cc|c|cccc|c|cc|c|cc}
0_{s_1}&&&&&&&&&&&&&\\\hline
&0&I_{s_2}&&&&&&&&&&&\\
&&0&&&&&&&&&&&\\\hline
&&&\ddots&&&&&&&&&&\\\hline
&&&&0&I_{s_\rho}&&&&&&&&\\
&&&&&0&\ddots&&&&&&&\\
&&&&&&\ddots&I_{s_\rho}&&&&&&\\
B_{\rho1}^{(\rho1)}&B_{\rho1}^{(\rho2)}&&\ldots&B_{\rho1}^{(\rho\rho)}&&&0&
W_{\rho,\rho+1}&K_{\rho,\rho+1}&&\ldots&K_{\rho k}&\\\hline
Z_{\rho+1,1}&Z_{\rho+1,2}&&\ldots&Z_{\rho+1,\rho}&&&&W_{\rho+1}&
Z_{\rho+1,\rho+1}&&\ldots&Z_{\rho+1,k}&\\\hline
&&&&&&&&&I_{s_{\rho+1}}&&&\\
&&&&&&&&&&I_{(\rho-1)s_{\rho+1}}&&\\\hline
&&&&&&&&&&&\ddots&&\\\hline
&&&&&&&&&&&&I_{(k-\rho)s_k}&\\
&&&&&&&&&&&&&I_{(\rho-1)s_k}
\rix,
$$
}
where $W_{\rho+1}$ is defined in (\ref{defw}),
\eqn
\nonumber
W_{\rho,\rho+1}&=&\mat{ccc}B_{\rho1}^{(\rho,\rho+1)}&\ldots&B_{\rho1}^{(\rho k)}\rix, \\
\nonumber
K_{\rho j}&=&\quad \mat{ccc}B_{\rho2}^{(\rho j)}&\ldots&
B_{\rho,j-\rho+1}^{(\rho j)}\rix,\quad j=\rho+1,\ldots,k,
\\
\label{blks1}
\wh I_j&=&\mat{c}0_{(s_{\rho+1}+\ldots+s_{j-1})\times s_{j}}
\\I_{s_{j}}\\0_{(s_{j+1}+\ldots+s_k)\times s_{j}}\rix,
\qquad j=\rho+1,\ldots,k,
\enn

\[
Z_{\rho+1,j}=\left\{\begin{array}{l}
\mat{c}B_{\rho+1, 1}^{(\rho+1, j)}\\\vdots\\B_{k1}^{(kj)}\rix,\qquad 1\le j\le\rho\\\\
\mat{ccc}
B_{\rho+1,2}^{(\rho+1, j)}&\ldots&B_{\rho+1,j-\rho+1}^{(\rho+1, j)}\\
\vdots&\vdots&\vdots\\
B_{k2}^{(k j)}&\ldots&
B_{k,j-\rho+1}^{(k j)}\rix,\qquad j=\rho+1\le j\le k.
\end{array}\right.
\]
With the assumption that $W_{\rho+1}$ is invertible, the blocks $Z_{\rho+1,j}$
to the left of $W_{\rho+1}$ in $\Pi_LV\Pi_R$ can be eliminated by postmultiplying 
\eq{defg}
G=
\mat{c|cc|c|cccc|c|c}
I_{s_1}&&&&&&&&&\\\hline
&I_{s_2}&&&&&&&&\\
&&I_{s_2}&&&&&&&\\\hline
&&&\ddots&&&&&&\\\hline
&&&&I_{s_\rho}&&&&&\\
&&&&&\ddots&&&&\\
&&&&&&I_{s_\rho}&&&\\
&&&&&&&I_{s_\rho}&&\\\hline
G_1^{(\rho)}&G_2^{(\rho)}&
&&G_\rho^{(\rho)}&&&&I_{\hat s_{\rho+1}}&\\\hline
&&&&&&&&&I_{\sum_{i=\rho+1}^k(i-1)s_{i}}
\rix,
\en
where 
\eq{hats}
\hat s_{\rho+1}=\sum_{i=\rho+1}^ks_i
\en
and 
\eq{blkg}
G_j^{(\rho)}=-W_{\rho+1}^{-1}Z_{\rho+1,j},\qquad j=1,\ldots,\rho.
\en

The elimination will not affect the matrix $\Pi_LU\Pi_R$. 
Then 
\eq{whuv}
\Pi_L(\mu U-V)\Pi_RG =
\mu\mat{ccc}I_{\sum_{i=1}^{\rho-1}is_i}&&\\&I_{\rho s_\rho}&\\&&U_{33}\rix-
\mat{ccc}V_{11}&0&0\\V_{21}&\Theta_\rho&V_{23}\\0&0&V_{33}\rix=:\mu\wh U-\wh V,
\en
where
\eqn
\nonumber
U_{33}&=&\mat{c|c|cc|c|cc}0_{\hat s_{\rho+1}}&\mat{cc}0&\wh I_{\rho+1}\rix&0&
\mat{cc}0&\wh I_{\rho+2}\rix&\ldots&0&\mat{cc}0&\wh I_k\rix
\\\hline
&\wh N_{\rho+1}^T&&&&&\\\hline
&&0_{s_{\rho+2}}&&&&\\
&&&\wh N_{\rho+2}^T&&&\\\hline
&&&&\ddots&&\\\hline
&&&&&0_{(k-\rho-1)s_k}&\\&&&&&&\wh N_k^T\rix,\\
\nonumber
V_{11}&=&\diag(N_1,N_2,\ldots,N_{\rho-1}),\\
\nonumber
V_{21}&=&\mat{c|cc|c|cc}0_{(\rho-1)s_\rho\times s_1}&0_{(\rho-1)s_\rho\times s_2}&
0_{(\rho-1)s_\rho\times s_2}&\ldots&0_{(\rho-1)s_\rho\times s_{\rho-1}}&
0_{(\rho-1)s_\rho\times s_{\rho-1}}\\
S_1&S_2&0_{s_\rho\times s_2}&\ldots&S_{\rho-1}&0_{s_\rho\times (\rho-2)s_{\rho-1}}\rix\\
\label{blks2}
\Theta_\rho&=&\mat{cccc}0&I_{s_\rho}&&\\&0&\ddots&\\
&&\ddots&I_{s_\rho}\\S_\rho&&&0\rix_{\rho\times \rho}\\
\nonumber
V_{23}&=&\mat{c|cc|c|cc}
0_{(\rho-1)s_\rho\times \hat s_{\rho+1}}
&0&0&
\ldots&
0&0\\
W_{\rho,\rho+1}&K_{\rho,\rho+1}&0_{\rho\times (\rho-1)s_{\rho+1}}&
\ldots&K_{\rho k}&0_{\rho\times (\rho-1)s_{k}}
\rix,\\
\nonumber
V_{33}&=&\mat{c|cc|c|cc}
W_{\rho+1}&Z_{\rho+1,\rho+1}&&&Z_{\rho+1,k}&\\\hline
&I_{s_{\rho+1}}&&&&\\
&&I_{(\rho-1) s_{\rho+1}}&&&\\\hline
&&&\ddots&&\\\hline
&&&&I_{(k-\rho)s_{k}}&\\
&&&&&I_{(\rho-1)s_k}\rix,\\
\nonumber
\wh N_i&=&\mat{cccc}0&I_{s_i}&&\\&0&\ddots&\\
&&\ddots&I_{s_i}\\&&&0\rix_{\rho\times \rho},\quad i=\rho+1,\ldots,k,\\
\nonumber
S_{i}&=&B_{\rho1}^{(\rho i)}+W_{\rho,\rho+1}G_i^{(\rho)},\qquad i=1,\ldots,\rho.
\enn
Note by using (\ref{defw}) and (\ref{blkg}),
\eq{factw}
W_\rho =\mat{cc}B_{\rho1}^{(\rho\rho)}&W_{\rho,\rho+1}\\Z_{\rho+1,\rho}&W_{\rho+1}\rix
=\mat{cc}S_\rho&W_{\rho,\rho+1}\\0&W_{\rho+1}\rix
\mat{cc}I&0\\-G_\rho^{(\rho)}&I\rix.
\en
Based on the block structure,
$\mu \wh U-\wh V$ is a regular pencil having the eigenvalues zero and infinity contained
in $\mu I-V_{11}$ and $\mu U_{33}-V_{33}$, respectively, and the other finite
eigenvalues contained in $\mu I-\Theta_\rho$.
Based on (\ref{factw}),
\eq{epen}
\gamma \mat{cc}I_{s_\rho}&0\\0&0\rix-W_\rho=
\left(\gamma\mat{cc}I_{s_\rho}&0\\0&0\rix-\mat{cc}S_\rho&W_{\rho,\rho+1}\\ 0&W_{\rho+1}\rix\right)\mat{cc}I&0\\-G_\rho^{(\rho)}&I\rix.
\en
So the $s_p$ finite eigenvalues of $\gamma \mat{cc}I_{s_\rho}&0\\0&0\rix-W_\rho$ are just the 
eigenvalues of $S_\rho$.
Let these $s_p$ eigenvalues be
$\gamma_1^{(\rho)},\ldots, \gamma_{s_\rho}^{(\rho)}$,
and for each $i$, $\mu_{i1}^{(\rho)},\ldots,\mu_{i\rho}^{(\rho)}$ be the $\rho$th roots of $\gamma_i^{(\rho)}$.
One can show that all these $\mu_{ij}^{(\rho)}$ are the eigenvalues of $\Theta_\rho$
and also eigenvalues of $\mu \wh U-\wh V$. Using the relations (\ref{eqmatrix}), (\ref{abuv}), and (\ref{whuv}), and the fact that $zE_U$ and $E_V(z)$ in (\ref{abuv}) are order $O(z)$, one can obtain
the eigenvalue perturbation results in Theorem~\ref{thm1}(a). Note that $P(t)=O(t)$, and then
$tD_{12}P(t)=O(t^2)$. Replacing $D_{11}$ by $D_{11}+D_{12}P(t)$ in $A_{11}+tD_{11}$ will introduce an additional 
term to each subblock of $E_V(z)$ (given in (\ref{abuv})) with order at least $z^\rho$, which
will not affect the eigenvalue perturbation results. For the same reason, we will identify $D_{11}$ with $D_{11}+D_{12}P(t)$ without repeating the fact in the processes of deriving our main results
in the next two sections.

 With the generic condition that all $W_1,\ldots,W_k$ 
are invertible, based on (\ref{factw}), $S_\rho$ is invertible and $\gamma_1^{(\rho)},\ldots, \gamma_{s_\rho}^{(\rho)}$ are nonzero. If in addition these $s_\rho$ eigenvalues are 
distinct, then all $\mu_{ij}^{(\rho)}$ are simple eigenvalues of $\Theta_\rho$.
Suppose $\varphi_i\ne 0$ is
an eigenvector of $S_\rho$ corresponding to $\gamma_i^{(\rho)}$. Then for any one of 
the $\rho$th roots of $\gamma_i^{(\rho)}$, say $\mu_{ij}^{(\rho)}$, one has
$$
\Theta_\rho \wt u_{ij}^{(\rho)}=\mu_{ij}^{(\rho)}\wt u_{ij}^{(\rho)},
\qquad
\wt u_{ij}^{(\rho)}=\mat{c}\varphi_i\\\mu_{ij}^{(\rho)}\varphi_i\\\vdots\\(\mu_{ij}^{(\rho)})^{\rho-1} \varphi_i\rix.
$$
By using (\ref{epen}), one has
\eq{invw}
\left(\gamma_i^{(\rho)}\mat{cc}I_{s_\rho}&0\\0&0\rix-W_\rho\right)
\mat{c}\varphi_i\\G_\rho^{(\rho)}\varphi_i\rix =0.
\en
Let
\eq{defup}
G_\rho^{(\rho)}\varphi_i=\mat{c}u_{\rho+1}\\\vdots\\u_{k}\rix,\quad
u_{j}\in\setC^{s_j},\quad j=\rho+1,\ldots,k.
\en
One can show that
$$
u_{ij}^{(\rho)}=\mat{c|c|c|c|cc|c|cc}0_{1\times s_1}&\ldots&0_{1\times (\rho-1)s_{\rho-1}}&
(\wt u_{i,j}^{(\rho)})^T&u_{\rho+1}^T&0_{1\times \rho s_{\rho+1}}&
\ldots&u_{k}^T&0_{1\times (k-1)s_k}\rix^T
$$ 
is an eigenvector of $\mu U-V$ corresponding to $\mu_{ij}^{(\rho)}$, i.e., $\mu_{ij}^{(\rho)}Uu_{ij}^{(\rho)}= Vu_{ij}^{(\rho)}.$ Applying  the standard perturbation theory
to (\ref{abuv}), one has
\[
\wh\mu_{ij}^{(\rho)}U(z)\wh u_{ij}^{(\rho)} = V(z)\wh u_{ij}^{(\rho)},\quad
\wh u_{ij}^{(\rho)}=u_{ij}^{(\rho)}+O(z),\qquad 
\wh \mu_{ij}^{(\rho)}=\mu_{ij}^{(\rho)}+O(z).
\]
Using (\ref{simad}), (\ref{eqmatrix}), and
(\ref{abuv}), the vector $x_{ij}^{(\rho)}=(\Xi+\Xi_cP(t)) R\wh u_{ij}^{(\rho)}$,
where $\Xi,\Xi_c,P(t)$ are given in (\ref{simad}) and $R$ is defined in (\ref{deflr}),
is an eigenvector of $A+tD$ corresponding to the eigenvalue 
$\lambda_0+t^{1/\rho} \mu_{ij}^{(\rho)}+O(t^{2/\rho})$. The eigenvector perturbation results in Theorem~\ref{thm1}(b) can be derived
essentially in this way by using the structure of $R$ and the fact that $P(t)=O(t)$. Note 
due to the power forms of the blocks in $R$, only $\varphi_i$ from $\wt u_{ij}^{(\rho)}$
and $u_{\rho+1},\ldots,u_k$ are involved in forming the constant term of $x_{ij}^{(\rho)}$,
which can be written as $\wt\Xi_\rho\mat{c}\varphi_i\\G_{\rho}^{(\rho)}\varphi_i\rix$.
Following (\ref{invw}), $\mat{c}\varphi_i\\G_{\rho}^{(\rho)}\varphi_i\rix$ is equivalent to 
$\nu_i^{(\rho)}$ given in Theorem~\ref{thm1}(b).
\begin{remark}\label{rem1}\rm
In \cite{Lid65}, the perturbation results were established by applying analytic function
theory.
\end{remark}

\section{Invariant subspace perturbation}\label{secmain}
We consider the perturbation of invariant subspaces when the matrix $A$ 
is perturbed to $A+tD$. More specifically, we 
study the perturbed invariant subspaces corresponding to the eigenvalues with
a $t^{1/\rho}$ fractional order perturbation added to $\lambda_0$ 
for a given $\rho\in\{1,\ldots,k\}$.
This can be viewed as a generalization of the eigenvector perturbation theory given in Theorem~\ref{thm1}. Due to (\ref{simad}) and the  arguments made in the previous section, we focus on the problem with the matrix $\lambda_0I+N+tD_{11}$,
where $N$ is given in (\ref{nila}) and $D_{11}$ is given in (\ref{blkb}).
To avoid the
complicated situation that the eigenvalues having the same fractional order perturbation are 
generated from Jordan blocks of sizes other than $\rho$, we require the
generic condition that all $W_1,\ldots,W_k$ are invertible.
We will adopt the notations that we used in the preview section. Following the relation
(\ref{abuv}) we begin with the invariant subspace of the matrix pencil $\mu  U(z)- V(z)$
corresponding to all the eigenvalues that are perturbed from those of $\Theta_\rho$, 
where $\Theta_\rho$ is given in (\ref{blks2}). Recall under the generic condition, 
both $S_\rho$ and $\Theta_\rho$ are invertible.


Let 
\eq{uvuv}
\mu\wh U(z)-\wh V(z)=\Pi_L(\mu U(z)-V(z))\Pi_RG
=\mu (\wh U+z\wh E_U)-(\wh V+\wh E_V(z)),
\en
where  $\Pi_L$ and $\Pi_R$ are the block permutations defined in the previous section,
$G$ is given in (\ref{defg}), $\mu \wh U-\wh V$ is defined in (\ref{whuv}), and 
\bstar
\mu z\wh E_U-\wh E_V(z)&:=&\Pi_L(\mu zE_U-E_V(z))\Pi_RG\\
&=&\mu z\mat{ccc}0&0&0\\0&0&0\\\wh E^U_{31}&\wh E^U_{32}&\wh E^U_{33}\rix
-\mat{ccc}\wh E_{11}^V(z)&\wh E_{12}^V(z)&\wh E_{13}^V(z)\\
\wh E_{21}^V(z)&\wh E_{22}^V(z)&\wh E_{23}^V(z)\\
\wh E_{31}^V(z)&\wh E_{32}^V(z)&\wh E_{33}^V(z)\rix.
\estar
Using the partitioning of the blocks $G_1^{(\rho)},\ldots,G_\rho^{(\rho)}$ defined in (\ref{blkg}):
\eq{partg}
\mat{cccc}G_1^{(\rho)}&G_2^{(\rho)}&\ldots&
G_\rho^{(\rho)}\rix=
\mat{cccc}G_{\rho+1,1}&G_{\rho+1,2}&\ldots&G_{\rho+1,\rho}\\
G_{\rho+2,1}&G_{\rho+2,2}&\ldots&G_{\rho+2,\rho}\\
\vdots&\vdots&\vdots&\vdots\\
G_{k1}&G_{k2}&\ldots&G_{k\rho}\rix,
\en
with $G_{ij}\in\setC^{s_i\times s_j}$ for $i=\rho+1,\ldots,k$, $j=1,\ldots,\rho$, one has
$$
\wh E^U_{31}=\mat{c|cc|c|cc}0&0&0&\ldots&0&0\\\hline
G_{\rho+1,1}&G_{\rho+1,2}&0_{s_{\rho+1}\times s_2}&\ldots&G_{\rho+1,\rho-1}
&0_{s_{\rho+1}\times (\rho-2)s_{\rho-1}}\\
0_{(\rho-1)s_{\rho+1}\times s_1}&0&0&\ldots&0&0\\\hline
G_{\rho+2,1}&G_{\rho+2,2}&0_{s_{\rho+2}\times s_2}&\ldots&G_{\rho+2,\rho-1}
&0_{s_{\rho+2}\times (\rho-2)s_{\rho-1}}\\
0_{\rho s_{\rho+2}\times s_1}&0&0&\ldots&0&0\\\hline
\vdots&\vdots&\vdots&\vdots&\vdots&\vdots\\\hline
G_{k1}&G_{k2}&0_{s_{k}\times s_2}&\ldots&G_{k,\rho-1}
&0_{s_{k}\times (\rho-2)s_{\rho-1}}\\
0_{(k-2)s_{k}\times s_1}&0&0&\ldots&0&0\\\hline
\rix,
$$
\bstar
\wh E^U_{32}&=&\mat{cc}0&0\\\hline
G_{\rho+1,\rho}&0_{s_{\rho+1}\times (\rho-1)s_\rho}\\
0_{(\rho-1)s_{\rho+1}\times s_\rho}&0\\\hline
G_{\rho+2,\rho}&0_{s_{\rho+2}\times (\rho-1)s_\rho}\\
0_{\rho s_{\rho+2}\times s_\rho}&0\\\hline
\vdots&\vdots\\\hline
G_{k\rho}&0_{s_{k}\times (\rho-1)s_\rho}\\
0_{(k-2)s_{k}\times s_\rho}&0\\\rix,\\\\
\wh E^U_{33}&=&\small\mat{c|cc|ccc|c|ccccc}
0_{\wh s_{\rho+1}}&&&&&&&&&&&\\\hline
\wh I_{{\rho+1}}^T&0_{s_{\rho+1}}&&&&&&&&&&\\
&&0_{(\rho-1) s_{\rho+1}}&&&&&&&&&\\\hline
\wh I_{{\rho+2}}^T&&&0&&&&&&&&\\
&&&I_{s_{\rho+2}}&0&&&&&&&\\
&&&&&0_{(\rho-1) s_{\rho+2}}&&&&&&\\\hline
\vdots&&&&&\ddots&&&&&\\\hline
\wh I_{k}^T&&&&&&&0&&&&\\
&&&&&&&I_{s_k}&0&&&\\
&&&&&&&&\ddots&\ddots&&\\
&&&&&&&&&I_{s_k}&0&\\
&&&&&&&&&&&0_{(\rho-1)s_k}\rix;
\estar
where $\wh I_{\rho+1}\ldots,\wh I_k$ are defined in (\ref{blks1}), and
\bstar
\small
\wh E^V_{11}(z)&=&\mat{c|cc|c|ccc}z^{\rho-1}\wh B_{11}^{(11)}&
z^{\rho-1}\wh B_{11}^{(12)}&O(z^\rho)&\ldots&z^{\rho-1}\wh B_{11}^{(1,\rho-1)}&\ldots&
O(z^{2\rho-3})\\\hline
z^{\rho-1}\wh B_{11}^{(21)}&
z^{\rho-1}\wh B_{11}^{(22)}&O(z^\rho)&\ldots&z^{\rho-1}\wh B_{11}^{(2,\rho-1)}&\ldots&
O(z^{2\rho-3})\\
z^{\rho-2}\wh B_{21}^{(21)}&
z^{\rho-2}\wh B_{21}^{(22)}&O(z^{\rho-1})&\ldots&z^{\rho-2}\wh B_{21}^{(2,\rho-1)}&\ldots&
O(z^{2\rho-4})\\\hline
\vdots&\vdots&\vdots&\vdots&\vdots&\vdots&\vdots\\\hline
z^{\rho-1}\wh B_{11}^{(\rho-1,1)}&
z^{\rho-1}\wh B_{11}^{(\rho-1,2)}&O(z^\rho)&\ldots
&z^{\rho-1}\wh B_{11}^{(\rho-1,\rho-1)}&\ldots&O(z^{2\rho-3})\\
\vdots&\vdots&\vdots&\vdots&\vdots&\vdots&\vdots\\
z\wh B_{\rho-1,1}^{(\rho-1,1)}&
z\wh B_{11}^{(\rho-1,2)}&O(z^2)&\ldots
&z\wh B_{\rho-1,1}^{(\rho-1,\rho-1)}&\ldots&O(z^{\rho-1})\rix,
\\\\
\wh E^V_{12}(z)&=&
\mat{cccc}z^{\rho-1}\wh B_{11}^{(1\rho)}&O(z^\rho)&\ldots&O(z^{2\rho-2})\\
\hline
z^{\rho-1}\wh B_{11}^{(2\rho)}&O(z^\rho)&\ldots&O(z^{2\rho-2})\\
z^{\rho-2}\wh B_{21}^{(2\rho)}&O(z^{\rho-1})&\ldots&O(z^{2\rho-3})\\\hline
\vdots&\vdots&\vdots&\vdots\\\hline
z^{\rho-1}\wh B_{11}^{(\rho-1,\rho)}&O(z^\rho)&\ldots&O(z^{2\rho-2})\\
\vdots&\vdots&\vdots&\vdots\\
z\wh B_{\rho-1,1}^{(\rho-1,\rho)}&O(z^2)&\ldots&O(z^{\rho})
\rix, 
\estar
\bstar
\wh E^V_{21}(z)&=&\mat{c|cc|c|cccc}
z^{\rho-1}\wh B_{11}^{(\rho1)}&z^{\rho-1}\wh B_{11}^{(\rho2)}&O(z^\rho)&\ldots&
z^{\rho-1}\wh B_{11}^{(\rho,\rho-1)}&O(z^{\rho})&\ldots&O(z^{2\rho-3})\\
\vdots&\vdots&\vdots&\vdots&\vdots&\vdots&\vdots\\
z\wh B_{\rho-1,1}^{(\rho1)}&z\wh B_{\rho-1,1}^{(\rho2)}&O(z^2)&\ldots&
z\wh B_{\rho-1,1}^{(\rho,\rho-1)}&O(z^{2})&\ldots&O(z^{\rho-1})\\
0&0&zB_{\rho2}^{(\rho2)}&\ldots&0&zB_{\rho2}^{(\rho,\rho-1)}&\ldots
&O(z^{\rho-2})\rix,\\\\ 
\wh E_{22}^V(z)&=&\mat{ccccc}z^{\rho-1}\wh B_{11}^{(\rho\rho)}&O(z^{\rho})&O(z^{\rho+1})&\ldots&
O(z^{2\rho-1})\\
\vdots&\vdots&\vdots&\vdots&\vdots\\
z\wh B_{\rho-1,1}^{(\rho\rho)}&O(z^2)&O(z^{3})&\ldots&
O(z^{\rho})\\
0&zB_{\rho2}^{(\rho\rho)}&O(z^{2})&\ldots&
O(z^{\rho-1})\rix,\\\\
\wh E_{31}^V(z)&=&{\small \mat{c|cc|c|cccc}
0&0&zB_{\rho+1,2}^{(\rho+1,2)}&\ldots&0&zB_{\rho+1,2}^{(\rho+1,\rho-1)}&\ldots&
O(z^{\rho-2})\\
\vdots&\vdots&\vdots&\vdots&\vdots&\vdots&\vdots&\vdots\\
0&0&zB_{k2}^{(k2)}&\ldots&0&zB_{k2}^{(k,\rho-1)}&\ldots&
O(z^{\rho-2})\\\hline
z^\rho \wh B_{11}^{(\rho+1,1)}&z^\rho\wh B_{11}^{(\rho+1,2)}&O(z^{\rho+1})&
\ldots&O(z^\rho)&O(z^{\rho+1})&\ldots&
O(z^{2\rho-2})\\
\vdots&\vdots&\vdots&\vdots&\vdots&\vdots&\vdots&\vdots\\
z\wh B_{\rho1}^{(\rho+1,1)}&z\wh B_{\rho1}^{(\rho+1,2)}&O(z^2)&\ldots&
z\wh B_{\rho1}^{(\rho+1,\rho-1)}&O(z^2)&\ldots&
O(z^{\rho-1})\\\hline
\vdots&\vdots&\vdots&\vdots&\vdots&\vdots&\vdots&\vdots\\\hline
z^\rho \wh B_{11}^{(k1)}&z^\rho\wh B_{11}^{(k2)}&O(z^{\rho+1})&
\ldots&z^\rho\wh B_{11}^{(k,\rho-1)}&O(z^{\rho+1} )&\ldots&
O(z^{2\rho-2})\\
\vdots&\vdots&\vdots&\vdots&\vdots&\vdots&\vdots&\vdots\\
z^\rho \wh B_{k-\rho,1}^{(k1)}&z^\rho\wh B_{k-\rho,1}^{(k2)}&O(z^{\rho+1})&
\ldots&z^\rho\wh B_{k-\rho,1}^{(k,\rho-1)}&O(z^{\rho+1} )&\ldots&
O(z^{2\rho-2})\\
\vdots&\vdots&\vdots&\vdots&\vdots&\vdots&\vdots&\vdots\\
z\wh B_{k-1,1}^{(k1)}&z\wh B_{k-1,1}^{(k2)}&O(z^2)&\ldots&
z\wh B_{k-1,1}^{(k,\rho-1)}&O(z^2)&\ldots&
O(z^{\rho-1})\rix,}
\estar
\bstar
\wh E_{32}^V(z)&=&\mat{ccccc}0&zB_{\rho+1,2}^{(\rho+1,\rho)}&O(z^2)&\ldots&
O(z^{\rho-1})\\
\vdots&\vdots&\vdots&\vdots\\
0&zB_{k2}^{(k\rho)}&O(z^2)&\ldots&
O(z^{\rho-1})\\\hline
z^\rho \wh B_{11}^{(\rho+1,\rho)}&O(z^{\rho+1})&O(z^{\rho+2})&\ldots&
O(z^{2\rho-1})\\
\vdots&\vdots&\vdots&\vdots&\vdots\\
z\wh B_{\rho1}^{(\rho+1,\rho)}&O(z^2)&O(z^3)&\ldots&
O(z^{\rho})\\\hline
\vdots&\vdots&\vdots&\vdots\\\hline
z^\rho\wh B_{11}^{(k\rho)}&O(z^{\rho+1} )&O(z^{\rho+2})&\ldots&
O(z^{2\rho-1})\\
\vdots&\vdots&\vdots&\vdots\\
z^\rho\wh B_{k-\rho,1}^{(k\rho)}&O(z^{\rho+1})&O(z^{\rho+2})&\ldots&
O(z^{2\rho-1})\\
\vdots&\vdots&\vdots&\vdots&\vdots\\
z \wh B_{k-1,1}^{(k\rho)}&O(z^2)&O(z^3)&\ldots&
O(z^{\rho})\rix,
\estar
{\tiny
\bstar
&&\wh E^V_{13}(z)=\\
&&\mat{ccc|ccc|c|ccccc}
O(z^{\rho-1})&\ldots&O(z^{\rho-1})&
O(z^{\rho-1})&\ldots&O(z^{2\rho-2})&\ldots&
O(z^{\rho-1})&\ldots &O(z^{\rho-1})&\ldots&O(z^{2\rho-2})\\\hline
O(z^{\rho-1})&\ldots&O(z^{\rho-1})&
O(z^{\rho-1})&\ldots&O(z^{2\rho-2})&\ldots&
O(z^{\rho-1})&\ldots &O(z^{\rho-1})&\ldots&O(z^{2\rho-2})\\
O(z^{\rho-2})&\ldots&O(z^{\rho-2})&
O(z^{\rho-2})&\ldots&O(z^{2\rho-3})&\ldots&
O(z^{\rho-2})&\ldots &O(z^{\rho-2})&\ldots&O(z^{2\rho-3})\\\hline
\vdots&\vdots&\vdots&\vdots&\vdots&\vdots&\vdots&\vdots&\vdots&\vdots&\vdots&\vdots\\\hline
O(z^{\rho-1})&\ldots&O(z^{\rho-1})&
O(z^{\rho-1})&\ldots&O(z^{2\rho-2})&\ldots&
O(z^{\rho-1})&\ldots &O(z^{\rho-1})&\ldots&O(z^{2\rho-2})\\
\vdots&\vdots&\vdots&\vdots&\vdots&\vdots&\vdots&\vdots&\vdots&\vdots&\vdots&\vdots\\
zB_{\rho-1,1}^{(\rho-1,\rho+1)}&\ldots&zB_{\rho-1,1}^{(\rho-1,k)}&
zB_{\rho-1,2}^{(\rho-1,\rho+1)}&\ldots&O(z^{\rho})&\ldots&
zB_{\rho-1,2}^{(\rho-1,k)}&\ldots &zB_{\rho-1,k-\rho+1}^{(\rho-1,k)}&\ldots&O(z^{\rho})
\rix,
\estar
\bstar
&&\wh E^V_{23}(z)=\\
&&\mat{ccc|cccc|c|ccccc}
O(z^{\rho-1})&\ldots&O(z^{\rho-1})&
O(z^{\rho-1})&O(z^\rho)&\ldots&O(z^{2\rho-2})&\ldots&
O(z^{\rho-1})&\ldots &O(z^\rho)&\ldots&O(z^{2\rho-2})\\
\vdots&\vdots&\vdots&\vdots&\vdots&\vdots&\vdots&\vdots&\vdots&\vdots&\vdots\\
zB_{\rho-1,1}^{(\rho,\rho+1)}&\ldots&zB_{\rho-1,1}^{(\rho k)}&
zB_{\rho-1,2}^{(\rho,\rho+1)}&O(z^2)&\ldots&O(z^{\rho})&\ldots&
zB_{\rho-1,2}^{(\rho k)}&\ldots &O(z^2)&\ldots&O(z^{\rho})\\
0&\ldots&0&0&zB_{\rho3}^{(\rho,\rho+1)}&\ldots&O(z^{\rho-1})&\ldots&
0&\ldots&zB_{\rho,k-\rho+2}^{(\rho k)}&\ldots&O(z^{\rho-1})
\rix,\\
&&\wh E^V_{33}(z)=\\
&&\mat{ccc|ccc|c|ccccc}
0&\ldots&0&0&\ldots&O(z^{\rho-1})&\ldots&0&\ldots&0&
\ldots&O(z^{\rho-1})\\
\vdots&\vdots&\vdots&\vdots&\vdots&
\vdots&\vdots&\vdots&\vdots&\vdots&\vdots&\vdots\\
0&\ldots&0&0&\ldots&O(z^{\rho-1})&\ldots&0&\ldots&0&
\ldots&O(z^{\rho-1})\\\hline
O(z^\rho)&\ldots&O(z^\rho)&O(z^\rho)&
\ldots&O(z^{2\rho-1})&\ldots&
O(z^\rho)&\ldots&
O(z^\rho )&
\ldots&O(z^{2\rho-1})\\
\vdots&\vdots&\vdots&\vdots&\vdots&
\vdots&\vdots&\vdots&\vdots&\vdots&\vdots&\vdots\\
z B_{\rho1}^{(\rho+1,\rho+1)}&\ldots&z B_{\rho1}^{(\rho+1,k)}
&z B_{\rho2}^{(\rho+1,\rho+1)}&
\ldots&O(z^{\rho})&\ldots&
z B_{\rho2}^{(\rho+1,k)}&\ldots&
z B_{\rho,k-\rho+1}^{(\rho+1,k)}&
\ldots&O(z^{\rho})\\ \hline
\vdots&\vdots&\vdots&\vdots&\vdots&
\vdots&\vdots&\vdots&\vdots&\vdots&\vdots&\vdots\\\hline
O(z^\rho)&\ldots&O(z^\rho)&O(z^\rho)&
\ldots&O(z^{2\rho-1})&\ldots&
O(z^\rho)&\ldots& O(z^\rho)&\ldots&O(z^{2\rho-1})\\
\vdots&\vdots&\vdots&\vdots&\vdots&
\vdots&\vdots&\vdots&\vdots&\vdots&\vdots&\vdots\\
O(z^\rho)&\ldots&O(z^\rho)&O(z^\rho)&\ldots&O(z^{2\rho-1})&\ldots&
O(z^\rho)&\ldots&O(z^\rho )&\ldots&O(z^{2\rho-1})\\
\vdots&\vdots&\vdots&\vdots&\vdots&
\vdots&\vdots&\vdots&\vdots&\vdots&\vdots&\vdots\\
z B_{k-1,1}^{(k,\rho+1)}&\ldots&z B_{k-1,1}^{(kk)}
&z B_{k-1,2}^{(k,\rho+1)}&
\ldots&O(z^{\rho})&\ldots&
z B_{k-1,2}^{(kk)}&\ldots&
z B_{k-1,k-\rho+1}^{(kk)}&
\ldots&O(z^{\rho})
\rix,
\estar
}
where $O(z^j)$ represents a block that is $z^j$ times a corresponding subblock of 
$D_{11}$, and 
$$
\wh B_{i1}^{(j\ell)}=B_{i1}^{(j\ell)}+\mat{ccc}B_{i1}^{(j,\rho+1)}&\ldots&B_{i1}^{(jk)}\rix 
G_\ell^{(\rho)},
\quad i=1,\ldots, i_{\max},\quad j=1,\ldots,k,\quad \ell=1,\ldots,\rho,
$$
with $i_{\max}=\left\{\begin{array}{ll}j,&j\le \rho-1\\ j-1,&j>\rho-1.\end{array}\right.$

Now 
$
\mu \wh U(z)-\wh V(z)
=
$
$$
\mu\mat{ccc}I&0&0\\0&I&0\\z\wh E_{31}^U&z\wh E_{32}^U&U_{33}+z\hat E_{33}^U\rix -
\mat{ccc}V_{11}+\wh E_{11}^V(z)&\wh E_{12}^V(z)&\wh E_{13}^V(z)
\\V_{21}+\wh E_{21}^V(z)&\Theta_\rho+\wh E_{22}^V(z)&V_{23}+\wh E_{23}^V(z)\\
\wh E_{31}^V(z)&\wh E_{32}^V(z)&V_{33}+\wh E_{33}^V(z)\rix.
$$
We need to determine a basis for the deflating subspace of this pencil corresponding to the
subpencil that is perturbed from $\mu I-\Theta_\rho$. For this we
 eliminate the block $z\wh E_{32}^U$ in the first matrix and
the blocks $\wh E_{12}^V(z)$ and $\wh E_{32}^V(z)$ in the second matrix
 simultaneously by 
multiplying the matrices
$$
\mat{ccc}I&Y_1&0\\0&I&0\\0&Y_2&I\rix,\quad
\mat{ccc}I&X_1&0\\0&I&0\\0&X_2&I\rix
$$
from the left and right  sides of the pencil, respectively. This requires the blocks
$X_1,X_2, Y_1,Y_2$ to satisfy the  equations
\bstar
&&X_1+Y_1=0,\\
&&Y_2+(U_{33}+z\wh E^U_{33})X_2+z\wh E_{31}^UX_1+z\wh E_{32}^U=0,\\
&&(V_{11}+\wh E_{11}^V(z))X_1+Y_1\wh\Theta_\rho
+\wh E_{12}^V(z)+\wh E_{13}^V(z)X_2=0,\\
&&
(V_{33}+\wh E_{33}^V(z))X_2+Y_2\wh \Theta_\rho+\wh E_{31}^V(z)X_1+\wh E_{32}^V(z)=0,
\estar
where
\eq{hatt}
\wh \Theta_\rho =\Theta_\rho+\wh E_{22}^V(z)+(V_{21}+\wh E_{21}^V(z))X_1
+(V_{23}+\wh E_{23}^V(z))X_2.
\en

By eliminating $Y_1$ and $Y_2$ one has the generalized Riccati equation 
\bstar
&&
\mat{cc}V_{11}+\wh E_{11}^V(z)&\wh E_{13}^V(z)\\
\wh E_{31}^V(z)-z\wh E_{32}^U(V_{21}+\wh E_{21}^V(z))&V_{33}+\wh E_{33}^V(z)-z\wh E_{32}^U(V_{23}+\wh E_{23}^V(z))\rix\mat{c}X_1\\X_2\rix\\
&&-\mat{cc}I&0\\z \wh E_{31}^U&U_{33}+z\wh E_{33}^U\rix\mat{c}X_1\\X_2\rix (\Theta_\rho+\wh E_{22}^V(z))+\mat{c}\wh E_{12}^V(z)\\\wh E_{32}^V(z)-z\wh E_{32}^U(\Theta_\rho+\wh E_{22}^V(z))\rix\\
&&-\mat{cc}I&0\\z \wh E_{31}^U&U_{33}+z\wh E_{33}^U\rix\mat{c}X_1\\X_2\rix
\mat{cc}V_{21}+\wh E_{21}^V(z)&V_{23}+\wh E_{23}^V(z)\rix\mat{c}X_1\\X_2\rix=0.
\estar
Since the pencil $\mu\mat{cc}I&0\\0&U_{33}\rix-\mat{cc}V_{11}&0\\0&V_{33}\rix$
has only the eigenvalues $0$ and $\infty$, and the pencil
$\mu I-\Theta_\rho$ has only nonzero eigenvalues (due to
the nonsingularity of $\Theta_\rho$), for $z$ sufficiently small, one has
$X_1,X_2=O(z)$ (\cite{SteS90}).  Rewrite  the Riccati equation as
\eqn
\label{eqo}
&&V_{11}X_1-X_1\wh\Theta_\rho+\wh E_{11}^V(z)X_1+\wh E_{12}^V(z)+\wh E_{13}^V(z)X_2=0,\\
\label{eqt}
&&V_{33}X_2-(U_{33}+z\wh E_{33}^U)X_2\wh\Theta_\rho+E_2=0,
\enn
where
\bstar
E_2&=&(\wh E^V_{31}(z)-z\wh E_{32}^U(V_{21}+\wh E_{21}^V(z)))X_1-
z\wh E_{31}^UX_1\wh \Theta_\rho\\
&&+(\wh E_{33}^V(z)-z\wh E_{32}^U\wh E_{23}^V(z))X_2+\wh E_{32}^V(z)-z\wh E_{32}^U(\Theta_\rho+\wh E_{22}^V(z)+V_{23}X_2).
\estar
Using the block structures, $\wh E_{32}^UV_{21}=0$ if $\rho\ge 2$ or it is void if $\rho=1$.
Similarly, $\wh E_{32}^UV_{23}=0$ if $\rho\ge 2$ or it is of order $O(1)$ if $\rho=1$. In the latter case
$\Theta_1=S_1$ and then $\Theta_1+V_{23}X_2=S_1+O(z)$.  In either case, we may drop the
 terms $-z\wh E^U_{32}V_{21}X_1$ and $-z\wh E_{32}^UV_{23}X_2$ in 
 $E_2$. 

From the block forms of  $\wh E_{11}^V(z)$, $\wh E_{12}^V(z)$ and $\wh E_{13}^V(z)$, and 
the fact $X_1,X_2=O(z)$, 
$$
\wh E_{11}^V(z)X_1+\wh E_{12}^V(z)+\wh E_{13}^V(z)X_2
=O(\mat{c|cc|c|ccc}z^{\rho-1}&
z^{\rho-1}& z^{\rho-2}&\ldots& z^{\rho-1}&\ldots&z\rix^T).
 $$
Partition 
$$
X_1=\mat{ccc}\wh X_{1}^T&\ldots&
\wh X_{\rho-1}^T\rix^T,\qquad X_2=\mat{cccc} \wh X_\rho^T&
\wh X_{\rho+1}^T&\ldots &\wh X_{k}^T\rix^T.
$$ 
By comparing the top blocks on both sides of (\ref{eqo}) and the using the block structure of
$V_{11}$, we have
$$
\wh X_{1}\wh \Theta_\rho=O(z^{\rho-1}),
$$
which implies $\wh X_{1}=O(z^{\rho-1})$. By comparing the next block on both sides of (\ref{eqo}), one has
$$
N_2\wh X_{2}-\wh X_{2}\wh \Theta_\rho
 =O\left(\mat{c}z^{\rho-1}\\z^{\rho-2}\rix\right),
$$
and from which we have $\wh X_{2}=O(z^{\rho-2})$. By induction, we have
$$
\wh X_{i}=O(z^{\rho-i}),\qquad i=1,\ldots,\rho-1.
$$

Moving on to the equation (\ref{eqt}), adding $-z\wh E_{32}^U\wh E_{21}^V(z)$ to
$\wh E_{31}^V(z)$ will not change the fractional orders of 
 the blocks in $\wh E_{31}^V(z)$.
So using $X_1=O(z)$,
\begin{align*}
&(\wh E^V_{31}(z)-z\wh E_{32}^U\wh E_{21}^V(z))X_1\\
&=O(\mat{ccc|cccc|c|ccccccc}
z^2&\ldots&z^2&z^{\rho+1}&z^{\rho}&\ldots&z^2&\ldots&z^{\rho+1}&\ldots&z^{\rho+1}&
z^\rho&z^{\rho-1}&\ldots&z^2\rix^T).
\end{align*}
Next,
$$-z\wh E_{31}^UX_1\wh\Theta_\rho
=O(\mat{ccc|cccc|c|ccccccc}
0&\ldots&0&z^2&0&\ldots&0&\ldots&z^2&0&\ldots&0&0&\ldots&0\rix^T).
$$
Using $X_2=O(z)$ and the fact that adding $-z\wh E_{32}^U\wh E_{23}^V(z)$ to
$\wh E_{33}^V(z)$ will not change the fractional orders of the blocks in 
$\wh E_{33}^V(z)$, one has
\begin{align*}
&(\wh E_{33}^V(z)-z\wh E_{32}^U\wh E_{23}^V(z))X_2\\
&=O(\mat{ccc|cccc|c|ccccccc}
z^2&\ldots&z^2&z^{\rho+1}&z^\rho&\ldots&z^2&\ldots&
z^{\rho+1}&\ldots&z^{\rho+1}&z^\rho&z^{\rho-1}&\ldots&z^2\rix^T).
\end{align*}
Clearly,
$$
\wh E_{32}^V(z)=O(\mat{ccc|cccc|c|ccccccc}
z&\ldots&z&z^\rho&z^{\rho-1}&\ldots&z&\ldots&z^\rho&\ldots&z^\rho&z^{\rho-1}&z^{\rho-2}&
\ldots&z\rix^T).
$$
Finally,
$$
-z\wh E_{32}^U(\Theta_\rho+\wh E_{22}^V(z))=
O(\mat{ccc|cccc|c|ccccccc}
0&\ldots&0&z&0&\ldots&0&\ldots&z&0&0&\ldots&0&\ldots&0\rix^T).
$$
Altogether, one has
$$
E_2=O(\mat{ccc|cccc|c|ccccccc}
z&\ldots&z&z&z^{\rho-1}&\ldots&z&\ldots&z&z^\rho&\ldots&z^\rho&z^{\rho-1}&\ldots&z\rix^T).
$$
Note
$
U_{33}+z\wh E_{33}^U=$
$$\small
\mat{c|cccc|ccccc|c|ccccccc}0_{\hat s_{\rho+1}}&0&\ldots&0&\wh I_{\rho+1}&0&0&
\ldots&0&\wh I_{\rho+2}&\ldots&0&\ldots&0&0&\ldots&0&\wh I_k
\\\hline
z\wh I_{\rho+1}^T&0&&&&&&&&&&&&&&&&\\
0&I_{s_{\rho+1}}&0&&&&&&&&&&&&&&&\\
\vdots&&\ddots&\ddots&&&&&&&&&&&&&&\\
0&&&I_{s_{\rho+1}}&0&&&&&&&&&&&&&\\\hline
z\wh I_{\rho+2}^T&&&&&0&&&&&&&&&&&&\\
0&&&&&zI_{s_{\rho+2}}&0&&&&&&&&&&&\\
0&&&&&&I_{s_{\rho+2}}&0&&&&&&&&&&\\
\vdots&&&&&&&\ddots&\ddots&&&&&&&&&\\
0&&&&&&&&I_{s_{\rho+2}}&0&&&&&&&&\\\hline
0&&&&&&&&&&\ddots&&&&&&&\\\hline
z\wh I_k^T&&&&&&&&&&&0&&&&&&\\
0&&&&&&&&&&&zI_{s_k}&0&&&&&\\
\vdots&&&&&&&&&&&&\ddots&\ddots&&&&\\
0&&&&&&&&&&&&&zI_{s_k}&0&&&\\
0&&&&&&&&&&&&&&I_{s_k}&0&&\\
\vdots&&&&&&&&&&&&&&&\ddots&\ddots&\\
0&&&&&&&&&&&&&&&&I_{s_k}&0
\rix.
$$

By comparing the blocks in (\ref{eqt}) from top to bottom and
using the block structures of $(U_{33}+z\wh E_{33}^U)$ and $V_{33}$, one has
$$
\wh X_\rho=O(\mat{ccc}z&\ldots&z\rix^T),\qquad
\wh X_{\rho+1}=O(\mat{ccc}z&\ldots&z\rix^T),
$$
$$
\wh X_{i}=\left\{\begin{array}{ll}O\mat{ccccccccc}z&\ldots&z^{i-\rho-1}&z^{i-\rho}&\ldots&z^{i-\rho}&
z^{i-\rho-1}&\ldots&z\rix^T),&\rho+2\le i\le 2\rho\\
O(\mat{ccccccccc}z&\ldots&z^{\rho-1}&z^\rho&\ldots&z^\rho&z^{\rho-1}&\ldots&
z\rix^T),&k\ge i>2\rho.\end{array}\right.
$$
By construction,
$$
\wh U(z)\mat{c}X_1\\I\\X_2\rix\wh \Theta_\rho
=\wh V(z)\mat{c}X_1\\I\\X_2\rix.
$$

We now connect this  to the corresponding invariant subspace of $A+tD$.
Using (\ref{abuv}) and (\ref{uvuv}) one has 
\eq{finaleq}
\Pi_LL(z\mu I-(N+z^\rho D_{11}))R\Pi_RG = \mu \wh U(z)-\wh V(z).
\en
Define 
$$
\wt X=R\Pi_RG\mat{c}X_1\\I\\ X_2\rix.
$$
One has
\eq{uvdef}
(N+z^\rho D_{11})\wt X=\wt X(z\wh \Theta_\rho).
\en
Using the structures of $R$, $\Pi_R$, and $G$, we have
\eq{wtx}
\wt X=\wt G+\mat{ccccccc}\wt X_{1}^T&\ldots&\wt X_{\rho-1}^T&\wt X_{\rho}^T&
\wt X_{\rho+1}^T&\ldots&\wt X_{k}^T\rix^T,
\en
where 
\bstar
\wt G&=&\mat{c}0_{s_1\times \rho s_\rho}\\\hline
\vdots\\\hline0_{(\rho-1)s_{\rho-1}\times \rho s_\rho}\\\hline
\mat{cc}I_{s_\rho}&0\\0&0_{(\rho-1)s_\rho\times (\rho-1) s_\rho}\rix\\\hline
\mat{cc}G_{\rho+1,\rho}&0\\0&0_{\rho s_{\rho+1}\times (\rho-1) s_\rho}\rix\\\hline
\vdots\\\hline \mat{cc}G_{k\rho}&0\\0&0_{(k-1) s_{k}\times (\rho-1) s_\rho}\rix\rix,\\
\wt X_{\rho}&=&\diag(0_{s_\rho},zI_{s_\rho},\ldots,z^{\rho-1}I_{s_\rho}),\\
\wt X_{i}
&=&\left\{\begin{array}{ll}
O(\mat{cccc}z^{\rho-i}&z^{\rho-i+1}&\ldots&z^{\rho-1}\rix^T),& 1\le i\le\rho-1\\
O(\mat{cccccccc}z& z&z^2&\ldots&z^{\rho-1}&z^\rho&\ldots&z^\rho\rix^T),
&\rho+1\le i\le k,\end{array}\right.
\estar
where $G_{\rho+1,\rho},\ldots,G_{k\rho}$ are the blocks of $G_\rho^{(\rho)}$ 
given in (\ref{partg}).
Then using $z=t^{1/\rho}$,  one has $(\lambda_0I+N+tD_{11})\wt X=\wt X(\lambda_0I+t^{1/\rho}\wh \Theta)$.
Following (\ref{simad}) (and identifying $D_{11}$ with $D_{11}+D_{12}P(t)$), one has
\[
(A+tD)X=X(\lambda_0I+t^{1/\rho}\wh\Theta_\rho),
\]
where
\eq{invab0} 
X=(\Xi+\Xi_cP(t))\wt X=\Xi\wt X+O(t).
\en
Hence, $X$ is a matrix of the invariant subspace of $A+tD$ corresponding to the 
eigenvalues of $\lambda_0I+t^{1/\rho}\wh \Theta_\rho$, i.e., all the eigenvalues perturbed
with a fractional order $t^{1/\rho}$. We summarize the results in the
following theorem.
\begin{theorem}\label{thm2.5}
Suppose $A$ and $D$ are transformed to (\ref{ta}) and (\ref{td}), where 
$A_{11}$ and $D_{11}$ are given in (\ref{nila}) and (\ref{blkb}),
respectively. Assume that $t>0$ is sufficiently small and $W_1,\ldots,W_\rho$ defined in (\ref{defw})
are invertible. Then for a given $\rho\in\{1,\ldots,k\}$ one has
\eq{invab}
(A+tD)X=X(\lambda_0I+t^{1/\rho}\wh \Theta_\rho),
\en
where $\wh \Theta_\rho$ is defined in (\ref{hatt}) satisfying
$\wh \Theta_\rho=\Theta_\rho+O(t^{1/\rho})$ with $\Theta_\rho$ given in (\ref{blks2}),
$$
X=\wt \Xi_\rho\mat{cc}I_{s_\rho}&0_{s_\rho\times (\rho-1)s_\rho}\\
G_\rho^{(\rho)}&0_{\hat s_{\rho+1}\times (\rho-1)s_\rho}\rix+\sum_{i=1}^k\Xi_i\wt X_{i}
+O(t),
$$
where $\hat s_{\rho+1}$ is defined in (\ref{hats}),
 $\wt \Xi_\rho$ is defined in (\ref{eigspace}), $G_\rho^{(\rho)}$ is given in (\ref{blkg}), 
$\Xi_1,\ldots,\Xi_k$ are the block columns of $\Xi$ defined in (\ref{bkxi}),
and
\bstar
\wt X_{\rho}&=&\diag(0_{s_\rho},t^{1/\rho}I_{s_\rho},\ldots,t^{1-1/\rho}I_{s_\rho}),\\
\wt X_{i}
&=&\left\{\begin{array}{ll}
O(\mat{cccc}t^{1-i/\rho}&t^{1-(i-1)/\rho}&\ldots&t^{1-1/\rho}\rix^T),& 1\le i\le\rho-1\\
O(\mat{cccccccc}t^{1/\rho}&t^{1/\rho}&t^{2/\rho}&\ldots&t^{1-1/\rho}&t&\ldots&t\rix^T),&\rho+1\le i\le k.\end{array}\right.
\estar
\end{theorem}
\proof
The expression of $X$ is simply from the formula  (\ref{invab0})
and the block forms of $\Xi$, $\wt\Xi$, and the partitioning of $G_\rho^{(\rho)}$ in (\ref{partg}).
The fractional order forms of $\wt X_{i}$ are obtained with the substitution $z=t^{1/\rho}$.
\eproof
The results of Theorem~\ref{thm2.5} reveal how  the eigenvectors and generalized
eigenvectors of $A$ corresponding to $\lambda_0$ are  involved in constructing $X$.

If the associated eigenvalues are perturbed with
different fractional orders, we may decouple the invariant subspace as a 
direct sum of several invariant subspaces with each corresponding to the eigenvalues with 
the same fractional order perturbation. We then apply the above results to each 
invariant subspace.

\medskip
Next, we consider the perturbation of an invariant subspace corresponding to part of the
eigenvalues of $\lambda_0I+t^{1/\rho}\wh \Theta_\rho$. For this we turn to consider 
an invariant subspace of $\wh\Theta_\rho$ that is given in (\ref{hatt}). 
Recall $\wh \Theta_\rho$ is perturbed from $\Theta_\rho$ given in (\ref{blks2}).
 Assume that $\wt \Phi_\rho\in\setC^{\rho s_\rho\times r}$
is a matrix of an invariant subspace of $\Theta_\rho$, i.e., $\rank \wt \Phi_\rho=r$ and
\eq{invphi}
\Theta_\rho\wt \Phi_\rho=\wt \Phi_\rho\Omega
\en
for some $\Omega\in\setC^{r\times r}$. Obviously, $\Lambda(\Omega)\subseteq\Lambda(\Theta_\rho)$.
We consider the case when the eigenvalues of $\Omega$ are different from the 
rest eigenvalues of
$\Theta_\rho$. From the block form of $\Theta_\rho$ one has
\eq{whp}
\wt\Phi_\rho=\mat{cccc}Q_1^T&(Q_1\Omega)^T&\ldots&(Q_1\Omega^{\rho-1})^T\rix^T
\en
for some matrix $Q_1\in\setC^{s_\rho\times r}$.  By comparing the 
bottom blocks on both sides of (\ref{invphi}) one has 
\eq{invs}
S_\rho Q_1= Q_1\Omega^{\rho}.
\en

On the other hand, suppose $\gamma^{(\rho)}_i$ is a single eigenvalue of $S_\rho$ with algebraic multiplicity 
$r_i$. Let $Q_1\in\setC^{s_\rho\times r_i}$ with $\rank Q_1=r_i$, and $S_{11}\in\setC^{r_i\times r_i}$ contain the only eigenvalue
$\gamma_i^{(\rho)}$ such that
$$
S_\rho Q_1=Q_1S_{11}.
$$
Since $S_\rho$ is invertible, so is $S_{11}$. Suppose that $\Omega$ is a $\rho$th
root of $S_{11}$ ($\Omega^{\rho}=S_{11}$) and contains a single eigenvalue
that is one of the $\rho$th roots of $\gamma_i^{(\rho)}$. Then the matrix 
$\wt \Phi_\rho$ constructed by the current $Q_1$ and $\Omega$ as defined in (\ref{whp}) 
has full column rank and satisfies (\ref{invphi}). 
Note the eigenstructure of the single eigenvalue of $\Omega$ 
is the same as that of $\gamma_i^{(\rho)}$ as an eigenvalue of $S_\rho$, in the sense
that they have the same partial multiplies (Jordan block sizes). 
The above relation is another way to
show that  $\Lambda(\Theta_\rho)$ is the collection of all the $\rho$th roots of the eigenvalues of $S_\rho$.

In general, $\Omega$ given in (\ref{invphi})  may have multiple eigenvalues. The assumption that the eigenvalues of
$\Omega$ are different from the rest eigenvalues of $\Theta_\rho$ implies that for
 each distinct eigenvalue of $\Omega$, a corresponding full column rank matrix as 
 $\wt \Phi_\rho$ can be constructed following the descriptions given above. The matrix $\wt\Phi_\rho$ 
 in (\ref{invphi}) can then be considered as formed by putting all such matrices side by side. 
Clearly, the resulting $\wt \Phi_\rho$ still has the form (\ref{whp}), and in this case $\Omega$ is a block 
diagonal matrix with each diagonal block containing a single eigenvalue.

The eigenvalues of $\Omega$ are different from the rest eigenvalues of $\Theta_\rho$,
and from (\ref{hatt}), $\wh\Theta_\rho=\Theta_\rho+O(z)$.
Following the standard invariant subspace perturbation theory, e.g., \cite{SteS90},
for $z$ sufficiently small,  there are matrix $\wh \Omega$ and full column rank matrix
$\wh\Phi_\rho$ such that
\eq{fstorder}
\wh\Theta_\rho\wh\Phi_\rho=\wh\Phi_\rho\wh \Omega,\qquad 
\wh \Phi_\rho=\wt\Phi_{\rho}+O(z), \quad\wh \Omega=\Omega+O(z).
\en
We have the following invariant subspace perturbation results.
\begin{theorem}\label{thm2}
Suppose $A$ and $D$ are transformed to (\ref{ta}) and (\ref{td}), where 
$A_{11}$ and $D_{11}$ are given in (\ref{nila}) and (\ref{blkb}),
respectively. Assume that $t>0$ is sufficiently small and $W_1,\ldots,W_\rho$ defined in (\ref{defw})
are invertible. Let $\Omega$ be a matrix that contains some eigenvalues of $\Theta_\rho$ that are different from the rest eigenvalues of $\Theta_\rho$. Then 
\eq{inab}
(A+tD)H=HC,
\en
where 
$$
C=\lambda_0I+t^{1/\rho}\Omega+O(t^{2/\rho})
$$
and
$$
H=\wt \Xi_\rho\mat{c}I_{s_\rho}\\G_\rho^{(\rho)}\rix Q_1+\sum_{i=1}^k
\Xi_iH_i+O(t)
$$
with $\wt\Xi_\rho$ given in (\ref{eigspace}), $\Xi_1,\ldots,\Xi_k$ defined in (\ref{bkxi}),
 $G_\rho^{(\rho)}$ given in (\ref{blkg}),  $Q_1$ given in (\ref{whp}),  
 and
\bstar
H_\rho&=&
\mat{cccc}0&(Q_1(t^{1/\rho}\Omega))^T&\ldots&
(Q_1(t^{1/\rho}\Omega)^{\rho-1})^T\rix^T+O(\mat{cccc}t^{1/\rho}&t^{2/\rho}&\ldots
&t\rix^T),\\
H_i&=&\left\{\begin{array}{ll} O(\mat{cccc}t^{1-i/\rho}&t^{1-(i-1)/\rho}&\ldots&
t^{1-1/\rho}\rix^T),&1\le i\le \rho-1\\
O(\mat{cccccccc}t^{1/\rho}& t^{1/\rho}&t^{2/\rho}&\ldots&t^{1-1/\rho}&t&\ldots&t\rix^T),
&\rho+1\le i\le k,
\end{array}\right.
\estar
\end{theorem}
\proof By postmultiplying $\wh \Phi_\rho$ on both sides of (\ref{invab}) and using 
(\ref{fstorder})  we have
$$
(A+tD)X\wh \Phi_\rho = X\wh \Phi_\rho (\lambda_0I+t^{1/\rho}\wh \Omega).
$$ 
Setting $C=\lambda_0I+t^{1/\rho}\wh \Omega$
and $H=X\wh\Phi_\rho$, we obtain (\ref{inab}). From (\ref{fstorder})
and using $z=t^{1/\rho}$, it is clear 
$$
C=\lambda_0I+t^{1/\rho}(\Omega+O(t^{1/\rho}))=
\lambda_0I+t^{1/\rho}\Omega+O(t^{2/\rho}),
$$
\[
H=\wt\Xi_\rho\mat{cc}I_{s_\rho}&0_{s_\rho\times (\rho-1)s_\rho}\\
 G_\rho^{(\rho)}&0_{\hat s_{\rho+1}\times (\rho-1)s_\rho}\rix(\wt\Phi_\rho+O(t^{1/\rho}))
 +\sum_{i=1}^k\Xi_i\wt X_i(\wt\Phi_\rho+O(t^{1/\rho})))+O(t).
 \]
Using the block structures of $\wt\Phi_\rho$ in (\ref{whp}),
\[
\wt\Xi_\rho\mat{cc}I_{s_\rho}&0_{s_\rho\times (\rho-1)s_\rho}\\
 G_\rho^{(\rho)}&0_{\hat s_{\rho+1}\times (\rho-1)s_\rho}\rix\wt\Phi_\rho =
 \wt\Xi_\rho\mat{c}I_{s_\rho}\\G_\rho^{(\rho)}\rix Q_1.
 \]
 Based on the block forms of $\wt \Xi_\rho$ and $\Xi_\rho,\ldots,\Xi_k$,
\bstar
 \wt\Xi_\rho\mat{cc}I_{s_\rho}&0_{s_\rho\times (\rho-1)s_\rho}\\
 G_\rho^{(\rho)}&0_{\hat s_{\rho+1}\times (\rho-1)s_\rho}\rix O(t^{1/\rho})
 &=&\wt\Xi_\rho\mat{ccc}O(t^{1/\rho})&\ldots& O(t^{1/\rho})\rix^T\\
&=&\sum_{i=\rho}^k\Xi_{i1}O(t^{1/\rho})
 =\sum_{i=\rho}^k\Xi_i\mat{cccc}O(t^{1/\rho})&0&\ldots&0\rix^T.
 \estar
 Then
  \[
 H_i=\left\{\begin{array}{ll}\wt X_i(\wt\Phi_\rho+O(t^{1/\rho})),&1\le i\le \rho-1\\
 \mat{cccc}O(t^{1/\rho})&0&\ldots&0\rix^T+\wt X_i(\wt\Phi_\rho+O(t^{1/\rho})),&\rho\le i\le k.
 \end{array}\right.
 \]
 The fractional orders of the blocks of $H_i$ can be derived based on the formulas, the forms of 
 $\wt X_i$ given in Theorem~\ref{thm2.5} and $\wt\Phi_\rho$ given in (\ref{whp}).
 \eproof
Note that following (\ref{epen}) and (\ref{invs}), one has
\[
W_\rho\mat{c}I_{s_\rho}\\G_\rho^{(\rho)}\rix Q_1
=\mat{cc}I_{s_\rho}&0\\0&0\rix \mat{c}I_{s_\rho}\\G_\rho^{(\rho)}\rix Q_1\Omega^\rho,
\]
which  generalizes (\ref{invw}). 
If $Q_1$ has full column rank, $\mat{c}I_{s_\rho}\\G_\rho^{(\rho)}\rix Q_1$ 
is a matrix of  the deflating subspace of the pencil $\mu \mat{cc}I_{s_\rho}&0\\0&0\rix-W_\rho$ 
corresponding to the eigenvalues of $\Omega^\rho$. 
However, $Q_1$ may be rank deficient.

\begin{remark}\rm $\mbox{}$
\begin{enumerate}
\item[(a)] If the eigenvalues of $\Omega$ in Theorem~\ref{thm2} have 
different fractional orders, perturbation results can be derived by
applying Theorem~\ref{thm2} to each set of the eigenvalues with the same fractional order.

\item[(b)] If  $\Lambda(\Omega)$ and $\Lambda(\Theta_\rho)\backslash\Lambda(\Omega)$
share common eigenvalues, we do not have (\ref{fstorder}).
In this situation, 
the problem becomes the original one about $A_{11}$ but with a reduced size.
So, we may apply Theorem~\ref{thm2} to each subproblem corresponding to
the eigenvalues of $\Omega$ with different fractional orders. 
If necessary, we may repeat
the process. Because the problem size is decreasing,  perturbation results can be 
derived eventually, if the generic condition always 
holds in the process.
\end{enumerate}
\end{remark}

We immediately have the following eigenvalue and eigenvector 
perturbation results. 
\begin{corollary}\label{cor1} 
Under the conditions of Theorem~\ref{thm2}, if $\gamma^{(\rho)}$ is a
simple eigenvalue of $S_\rho$
and $S_\rho \varphi=\gamma^{(\rho)} \varphi$ with $\varphi \ne 0$. Suppose $\mu^{(\rho)}$ satisfies $(\mu^{(\rho)})^\rho=\gamma^{(\rho)}$. Then $A+tD$ has an eigenvalue 
$$
\lambda =\lambda_0+t^{1/\rho}\mu^{(\rho)}+O(t^{2/\rho})
$$
and a corresponding eigenvector
$$
h=\wt\Xi\mat{c}\varphi\\G_\rho^{(\rho)}\varphi\rix+
\sum_{i=1}^k\Xi_ih_i+O(t),
$$
where $\wt\Xi_\rho$ is given in (\ref{eigspace}), $\Xi_1,\ldots,\Xi_k$ are defined in (\ref{bkxi}),
 $G_\rho^{(\rho)}$ is given in (\ref{blkg}), and
\bstar
h_\rho&=&\mat{cccc}0&(t^{1/\rho}\mu^{(\rho)})\varphi^T &\ldots&(t^{1/\rho}
\mu^{(\rho)})^{\rho-1}\varphi^T \rix^T
+O(\mat{cccc}t^{1/\rho}&t^{2/\rho}&\ldots&t\rix^T),\\
h_i
&=&\left\{\begin{array}{ll}
O(\mat{cccc}t^{1-i/\rho}&t^{1-(i-1)/\rho}&\ldots&t^{1-1/\rho}\rix^T),&1\le i\le\rho-1
\\
O(\mat{cccccccc}t^{1/\rho}& t^{1/\rho}&t^{2/\rho}&\ldots&t^{1-1/\rho}&t
&\ldots&t\rix^T),&\rho+1\le i\le k.\end{array}\right.
\estar
\end{corollary}
\proof
It is a direct consequence of Theorem~\ref{thm2} with $Q_1=\varphi$ and $\Omega=\mu^{(\rho)}$. 
\eproof
Based on (\ref{invw}), the vector $\mat{c}\varphi\\G_\rho^{(\rho)}\varphi\rix$ is an eigenvector
of the pencil (\ref{pen}) corresponding to the eigenvalue $\gamma^{(\rho)}$. Therefore,
by taking $\nu_i^{(\rho)}=\mat{c}\varphi\\G_\rho^{(\rho)}\varphi\rix$,
the constant term of $h$ is the same as the constant term of $x_{ij}^{(\rho)}$ in Theorem~\ref{thm1}(b). On the other hand, $h$ provides
more detailed perturbation information than $x_{ij}^{(\rho)}$ in the sum term.
Recall that under the assumptions of Corollary~\ref{cor1}, $\mu^{(\rho)}$ is a simple
eigenvalue of $\Theta_\rho$. Therefore, both $\lambda$ and 
$h$ have a series form as that in Theorem~\ref{thm1}. Note the series form 
 in Theorem~\ref{thm1}(b) exists for the eigenvectors associated with any
 simple eigenvalue $\gamma^{(\rho)}$ of $S_\rho$ or the pencil (\ref{pen}). The assumption 
 that all the finite eigenvalues of (\ref{pen}) are distinct is not necessary, unless one wants
to have a series form for the eigenvectors of all the order $t^{1/\rho}$ perturbed eigenvalues.

\medskip
We also have the  following results for the extreme cases.
\begin{corollary}\label{cor2} Suppose $\rho=1$ in Theorem~\ref{thm2}. Then one has
$$
(A+tD)H=HC,
$$
where $C=\lambda_0I+t\Omega+O(t^2)$ and $\Omega$ is defined in Theorem~\ref{thm2}, and
$$
H = \wt\Xi_1\mat{c}I_{s_1}\\G_1^{(1)}\rix Q_1+ O(t),
$$
with $\wt\Xi_1$ given in (\ref{eigspace}), $G_1^{(1)}$ given in (\ref{blkg}), and
$Q_1$ given in (\ref{whp}),  
\end{corollary}
\begin{corollary}\label{cor3} Suppose $\rho=k$ in Theorem~\ref{thm2}. Then one has
$$
(A+tD)H=HC,
$$
where $C=\lambda_0I+t^{1/k}\Omega+O(t^{2/k})$ and $\Omega$ is defined in 
Theorem~\ref{thm2},
$$
H = \Xi_{k1}Q_1+\sum_{i=1}^k\Xi_iH_i+ O(t)
$$
with $\Xi_{k1}(=\wt\Xi_k)$ and $\Xi_1,\ldots,\Xi_k$ given in (\ref{bkxi}),  $Q_1$
given in (\ref{whp}), and
$$
H_i=O(\mat{cccc}t^{1-i/k}&t^{1-(i-1)/k}&\ldots&t^{1-1/k}\rix^T),\quad 
1\le i\le k-1,
$$
$$
H_k=\mat{cccc}0&(Q_1(t^{1/k}\Omega))^T&\ldots&
(Q_1(t^{1/k}\Omega)^{k-1})^T\rix^T+O(\mat{cccc}t^{1/k}&t^{2/k}&\ldots&t\rix^T).
$$
\end{corollary}

%
\section{First fractional order invariant subspace perturbation}\label{secfst}
In Theorem~\ref{thm2}, we have 
\eq{ordh}
C=\lambda_0I+t^{1/\rho}\Omega+t^{2/\rho}\Omega_1+O(t^{3/\rho}),\quad
H=\wt \Xi_\rho \mat{c}I_{s_\rho}\\G_\rho\rix Q_1+t^{1/\rho}H_1+O(t^{2/\rho}).
\en
In this section,
we will derive
explicit expressions for $\Omega_1$ and $H_1$ by working on the pencil $\mu U(z)-V(z)$
given in (\ref{abuv}). Observe that the blocks in $V(z)$ with
an order higher than $z$ can be ignored in this case. So, we start with the pencil
$$
\mu (U+zE_U)-(V+zV_1),
$$
where $U,V,E_U$ are defined in (\ref{abuv}), and $V_1=$
{\tiny
$$
\mat{c|c|cccc|cccc|c|cccccc}
0_{s_1}&&&&&&&&&&&&&&&&\\\hline
&\ddots&&&&&&&&&&&&&&&\\\hline
&&0&&&&&&&&&&&&&&\\
&&&\ddots&&&&&&&&&&&&&\\
&&&&\ddots&&&&&&&&&&&&\\
B_{\rho-1,1}^{(\rho-1,1)}&\ldots&B_{\rho-1,1}^{(\rho-1,\rho-1)}&&&0&
B_{\rho-1,1}^{(\rho-1,\rho)}&&&&\ldots&B_{\rho-1,1}^{(\rho-1,k)}&\ldots&B_{\rho-1,k-\rho+1}^{(\rho-1,k)}&&&\\\hline
&&&&&&0&&&&&&&&&&\\
&&&&&&&\ddots&&&&&&&&&\\
B_{\rho-1,1}^{(\rho1)}&\ldots&B_{\rho-1,1}^{(\rho,\rho-1)}&&&&B_{\rho-1,1}^{(\rho\rho)}&&0&
&\ldots&B_{\rho-1,1}^{(\rho k)}&\ldots&B_{\rho-1,k-\rho+1}^{(\rho k)}&&&\\
0&\ldots&0&B_{\rho2}^{(\rho,\rho-1)}&&&0&B_{\rho2}^{(\rho\rho)}&&0&\ldots&
0&\ldots&0&B_{\rho,k-\rho+2}^{(\rho k)}&&\\\hline
\vdots&&&\vdots&&&&\vdots&&&\ddots&&&\vdots&&&\\\hline
&&&&&&&&&&&0&&&&&\\
&&&&&&&&&&&&0&&&&\\
&&&&&&&&&&&&&\ddots&&&\\
&&&&&&&&&&&&&&0&&\\
&&&&&&&&&&&&&&&0&\\
B_{k-1,1}^{(k1)}&\ldots&B_{k-1,1}^{(k,\rho-1)}&&&&B_{k-1,1}^{(k\rho)}&&&&\ldots&B_{k-1,1}^{(kk)}&
\ldots&B_{k-1,k-\rho+1}^{(kk)}&&\ddots&\\
0&\ldots&0&B_{k2}^{(k,\rho-1)}&&&0&B_{k2}^{(k\rho)}&&&\ldots&0&\ldots&0&B_{k,k-\rho+2}^{(kk)}&&0
\rix.
$$
}
We then apply the same techniques developed in the previous section to this pencil. 
We will use essentially the same notations used in the previous section,
although $E_V(z)$ in $V(z)$ is replaced by $zV_1$ here. 
We first transform the pencil to 
$$
\mu(\wh U+z\wh E_U)-(\wh V+z\wh V_1)=\Pi_L(\mu(U+zE_U)-(V+zV_1))\Pi_RG,
$$
where $\wh U, \wh V, \wh E_U$ are given  in the previous section,
$$
\wh V_1=\mat{ccc} E_{11}&E_{12}&E_{13}\\E_{21}&E_{22}&E_{23}\\E_{31}&E_{32}&E_{33}\rix
$$
and 
\bstar
E_{12}&=&\mat{cccc}0&0&\ldots&0\\\hline
\vdots&\vdots&\vdots&\vdots\\\hline
0&0&\ldots&0\\
\vdots&\vdots&\vdots&\vdots\\
0&0&\ldots&0\\\hline
0&0&\ldots&0\\
\vdots&\vdots&\vdots&\vdots\\
0&0&\ldots&0\\
\wh B_{\rho-1,1}^{(\rho-1,\rho)}&0&\ldots&0\rix,
\qquad E_{22}=\mat{ccccc}
0&0&0&\ldots&0\\
\vdots&\vdots&\vdots&\vdots&\vdots\\
0&0&0&\ldots&0\\
\wh B_{\rho-1,1}^{(\rho\rho)}&0&0&\ldots&0\\
0& B_{\rho2}^{(\rho\rho)}&0&\ldots&0\rix,
\estar
\bstar
&&E_{32}=\mat{ccccc}
0&B_{\rho+1,2}^{(\rho+1,\rho)}&0&\ldots&0\\
\vdots&\vdots&\vdots&\vdots&\vdots\\
0&B_{k2}^{(k\rho)}&0&\ldots&0\\\hline
0&0&0&\ldots&0\\
\vdots&\vdots&\vdots&\vdots&\vdots\\
0&0&0&\ldots&0\\
\wh B_{\rho1}^{(\rho+1,\rho)}&0&0&\ldots&0\\\hline
\vdots&\vdots&\vdots&\vdots&\vdots\\\hline
0&0&0&\ldots&0\\
\vdots&\vdots&\vdots&\vdots&\vdots\\
0&0&0&\ldots&0\\
\wh B_{k-1,1}^{(k\rho)}&0&0&\ldots&0\rix,
\estar
where
\bstar
\wh B_{\rho-1,1}^{(\rho-1,\rho)}&=& B_{\rho-1,1}^{(\rho-1,\rho)}+\mat{ccc}B_{\rho-1,1}^{(\rho-1,\rho+1)}&\ldots&
B_{\rho-1,1}^{(\rho-1,k)}\rix G_\rho^{(\rho)},\\
\wh B_{j-1,1}^{(j\rho)}&=& B_{j-1,1}^{(j\rho)}+\mat{ccc}B_{j-1,1}^{(j,\rho+1)}&\ldots&
B_{j-1,1}^{(jk)}\rix G_\rho^{(\rho)},\quad j=\rho,\ldots,k.
\estar
We do not include the explicit forms for other blocks of $\wh V_1$, since they are not
important in the following analysis. 

Now, by considering the first order only, the equations (\ref{eqo}) and (\ref{eqt}) become
two generalized Sylvester equations
\bstar
&&V_{11}X_1-X_1\Theta_\rho+zE_{12}=0,\\
&&V_{33}X_2-U_{33}X_2\Theta_\rho+z(E_{32}-\wh E_{32}^U\Theta_\rho)=0.
\estar
The first equation gives
$$
X_1=\mat{cccc}\wh X_{1}^T&\ldots&\wh X_{\rho-2}^T&\wh X_{\rho-1}^T\rix^T,
$$
where 
$$
\wh X_{i}=0_{is_i\times \rho s_\rho},\quad 1\le i\le \rho-2, \quad 
\wh X_{\rho-1}=z\mat{ccccc}0&\wh B_{\rho-1,1}^{(\rho-1,\rho)}S_\rho^{-1}&0&\ldots&0\\
\vdots&\vdots&\ddots&\vdots&\vdots\\
0&\ldots&0&\wh B_{\rho-1,1}^{(\rho-1,\rho)}S_\rho^{-1}&0\\
0&\ldots&0&0&\wh B_{\rho-1,1}^{(\rho-1,\rho)}S_\rho^{-1}\rix.
$$
In the second Sylvester equation
$$
E_{32}-\wh E_{32}^U\Theta_\rho =
\mat{ccccc}
0&B_{\rho+1,2}^{(\rho+1,\rho)}&0&\ldots&0\\
\vdots&\vdots&\vdots&\vdots&\vdots\\
0&B_{k2}^{(k\rho)}&0&\ldots&0\\\hline
0&-G_{\rho+1,\rho}&0&\ldots&0\\
\vdots&\vdots&\vdots&\vdots&\vdots\\
0&0&0&\ldots&0\\
\wh B_{\rho1}^{(\rho+1,\rho)}&0&0&\ldots&0\\\hline
\vdots&\vdots&\vdots&\vdots&\vdots\\\hline
0&-G_{k\rho}&0&\ldots&0\\
\vdots&\vdots&\vdots&\vdots&\vdots\\
0&0&0&\ldots&0\\
\wh B_{k-1,1}^{(k\rho)}&0&0&\ldots&0\rix.
$$
Again, let $X_2=\mat{cccc}\wh X_\rho^T&\wh X_{\rho+1}^T&\ldots&\wh X_{k}^T\rix^T$. 
From  the last block of the equation, we have
$$
\wh X_{k}-\mat{cc}0_{k-\rho-1}&0\\0&\wh N_k^T\rix
\wh X_{k}\Theta_\rho
+z\mat{ccccc}0&-G_{k\rho}&0&\ldots&0\\
0&0&0&\ldots&0\\
\vdots&\vdots&\vdots&\vdots&\vdots\\
0&0&0&\ldots&0\\
\wh B_{k-1,1}^{(k\rho)}&0&0&\ldots&0\rix=0,
$$
which gives
$$
\wh X_{k}=z\mat{ccccc}0&G_{k\rho}&0&\ldots&0\\
0&0&0&\ldots&0\\
\vdots&\vdots&\vdots&\vdots&\vdots\\
0&0&0&\ldots&0\\
-\wh B_{k-1,1}^{(k\rho)}&0&0&\ldots&0\rix.
$$
Similarly, we can derive
$$
\wh X_{i}=z\mat{ccccc}0&G_{i\rho}&0&\ldots&0\\
0&0&0&\ldots&0\\
\vdots&\vdots&\vdots&\vdots&\vdots\\
0&0&0&\ldots&0\\
-\wh B_{i-1,1}^{(i\rho)}&0&0&\ldots&0\rix,\quad \rho+2\le i\le k.
$$
The block $\wh X_{\rho+1}$ satisfies
$$
\wh X_{\rho+1}-\wh N_{\rho+1}^T \wh X_{\rho+1}\Theta_\rho
+z\mat{ccccc}0&-G_{\rho+1,\rho}&0&\ldots&0\\
0&0&0&\ldots&0\\
\vdots&\vdots&\vdots&\vdots&\vdots\\
0&0&0&\ldots&0\\
\wh B_{\rho1}^{(\rho+1,\rho)}&0&0&\ldots&0\rix=0,
$$
which gives
$$
\wh X_{\rho+1}=z\mat{ccccc}0&G_{\rho+1,\rho}&0&\ldots&0\\
0&0&G_{\rho+1,\rho}&\ldots&0\\
\vdots&\vdots&\vdots&\ddots&\vdots\\
0&0&0&\ldots&G_{\rho+1,\rho}\\
G_{\rho+1,\rho}S_\rho-\wh B_{\rho1}^{(\rho+1,\rho)}&0&0&\ldots&0\rix.
$$
Finally, $\wh X_\rho$ satisfies
$$
W_{\rho+1}\wh X_\rho=z\mat{ccccc}0&\wt C&0&\ldots&0\rix,
$$
where
$$
\wt C=\mat{c}G_{\rho+1,\rho}S_\rho\\0\\\vdots\\0\rix
-\mat{c|c|ccc}\wh B_{\rho1}^{(\rho+1,\rho)}&B_{\rho+1,2}^{(\rho+1,\rho)}&
B_{\rho+1,2}^{(\rho+1,\rho+1)}&\ldots&B_{\rho+1,2}^{(\rho+1,k)}\\
\wh B_{\rho+1,1}^{(\rho+2,\rho)}&B_{\rho+2,2}^{(\rho+2,\rho)}&
B_{\rho+2,2}^{(\rho+2,\rho+1)}&\ldots&B_{\rho+2,2}^{(\rho+2,k)}\\
\vdots&\vdots&\vdots&\vdots&\vdots\\
\wh B_{k-1,1}^{(k\rho)}&B_{k2}^{(k\rho)}&
B_{k2}^{(k,\rho+1)}&\ldots&B_{k2}^{(kk)}\rix
\mat{c}I\\\hline I\\\hline G_\rho^{(\rho)}\rix.
$$
Then
$$
\wh X_\rho
=z\mat{ccccc}0&\wh C&0&\ldots&0\rix,\quad \wh C=W_{\rho+1}^{-1}\wt C.
$$
Therefore,
$$
F:=\Pi_RG\mat{c}X_1\\I\\X_2\rix
=\wt G+\mat{ccccccccc}0_{s_1\times \rho s_\rho}^T&\ldots&0_{(\rho-2)s_{\rho-2}\time  \rho s_\rho}^T&F_{\rho-1}^T&F_{\rho}^T&F_{\rho+1}^T&F_{\rho+2}^T&\ldots& F_{k}^T\rix^T,
 $$
 where $\wt G$ is defined in (\ref{wtx}) and
 $$
 F_{\rho-1}=
 z\mat{ccccc}0&\wh B_{\rho-1,1}^{(\rho-1,\rho)}S_\rho^{-1}&0&\ldots&0\\
\vdots&\vdots&\ddots&\vdots&\vdots\\
0&\ldots&0&\wh B_{\rho-1,1}^{(\rho-1,\rho)}S_\rho^{-1}&0\\
0&\ldots&0&0&\wh B_{\rho-1,1}^{(\rho-1,\rho)}S_\rho^{-1}\rix,\quad
F_{\rho}=\mat{cc}0_{s_\rho}&0\\0&I_{(\rho-1)s_\rho}\rix,
$$
 $$
 F_{\rho+1}=z\mat{ccccc}
0&C_{\rho+1}&0&\ldots&0\\
 0&G_{\rho+1,\rho}&0&\ldots&0\\
0&0&G_{\rho+1,\rho}&\ldots&0\\
\vdots&\vdots&\vdots&\ddots&\vdots\\
0&0&0&\ldots&zG_{\rho+1,\rho}\\
G_{\rho+1,\rho}S_\rho-\wh B_{\rho1}^{(\rho+1,\rho)}&0&0&\ldots&0\rix,
$$
and
$$
F_{i}=z\mat{ccccc}0&C_{i}&0&\ldots&0\\
0&G_{i\rho}&0&\ldots&0\\
0&0&0&\ldots&0\\
\vdots&\vdots&\vdots&\vdots&\vdots\\
0&0&0&\ldots&0\\
-\wh B_{i-1,1}^{(i\rho)}&0&0&\ldots&0\rix,\quad \rho+2\le i\le k
$$
with
\eq{defr}
C=\wh C+G_{\rho-1}^{(\rho)}\wh B_{\rho-1,1}^{(\rho-1,\rho)}S_\rho^{-1}
=:\mat{c}C_{\rho+1}\\\vdots\\C_k\rix.
\en

From (\ref{uvdef}), we have
$$
(U+zE_U)F\wh \Theta_\rho = (V+zV_1)F+O(z^2),
$$
where  based on  (\ref{hatt}),
$$
\wh \Theta_\rho = \Theta_\rho+\Delta +O(z^2), 
$$
and
$$
\Delta =zE_{22}+V_{21}X_1+V_{23}X_2=z\mat{ccccc}0&0&0&\ldots&0\\
\vdots&\vdots&\vdots&\vdots&\vdots\\
0&0&0&\ldots&0\\
\wh B_{\rho-1,1}^{(\rho\rho)}&0&0&\ldots&0\\
0&\wh B_{\rho2}^{(\rho\rho)}&0&\ldots&0\rix,
$$
where
$$
\wh B_{\rho2}^{(\rho\rho)}=B_{\rho2}^{(\rho\rho)}+B_{\rho1}^{(\rho,\rho-1)}B_{\rho-1,1}^{(\rho-1,\rho)}S_\rho^{-1}+W_{\rho,\rho+1}C
+\mat{ccc}B_{\rho2}^{(\rho,\rho+1)}&\ldots&B_{\rho2}^{(\rho k)}\rix G_\rho^{(\rho)},
$$
and $C$ is defined in (\ref{defr}).
Let $\wt \Phi_\rho$ be defined in (\ref{whp}) and let $\wt \Phi_c$ be a  matrix
of  the invariant subspace corresponding to the eigenvalues of $\Theta_\rho$
other than those of $\Omega$. That is,
$$
\Theta_\rho \mat{cc}\wt \Phi_\rho&\wt \Phi_c\rix
=\mat{cc}\wt \Phi_\rho&\wt\Phi_c \rix\mat{cc}\Omega&0\\0&\Omega_c\rix
$$
and $\Lambda(\Omega)\cap\Lambda(\Omega_c)=\emptyset.$
Similar to $\wt \Phi_\rho$, $\wt \Phi_c$ has an expression
\eq{q2}
\wt \Phi_c=\mat{cccc}Q_2^T&(Q_2\Omega_c)^T&\ldots&(Q_2\Omega_c^{\rho-1})^T\rix^T
\en
for some matrix $Q_2$ that satisfies $S_\rho Q_2=Q_2 \Omega_c^\rho$.
Clearly, $\mat{cc}\wt \Phi_\rho&\wt\Phi_c\rix$ is invertible.
Let
$$
\mat{cc}\wt \Phi_\rho&\wt\Phi_c\rix^{-1}=\mat{c}\wt\Psi_\rho\\\wt \Psi_c\rix
$$
where $\wt \Psi_\rho$ is $r\times \rho s_\rho$. Analogously, one
has
\eq{moq}
\wt\Psi_\rho=
M\mat{cccc}\Omega^{\rho-1}\wt Q_1&\ldots&\Omega \wt Q_1&\wt Q_1\rix,\quad
\wt\Psi_c=M_c\mat{cccc}\Omega_c^{\rho-1}\wt Q_2&\ldots&\Omega_c \wt Q_2&\wt Q_2\rix,
\en
for some $\wt Q_1$, $\wt Q_2$ that satisfy
$$
\wt Q_1S_\rho=\Omega^\rho \wt Q_1,\quad
\wt Q_2S_\rho=\Omega_c^\rho\wt Q_2,
$$
and
$$
M=\left(\sum_{j=0}^{\rho-1}\Omega^{\rho-1-j}\wt Q_1Q_1\Omega^j\right)^{-1},\quad
M_c=\left(\sum_{j=0}^{\rho-1}\Omega_c^{\rho-1-j}\wt Q_2Q_2\Omega_c^j\right)^{-1}.
$$
Then
\[
\mat{c}\wt\Psi_\rho\\\wt\Psi_c\rix\wh \Theta_\rho\mat{cc}\wt\Phi_\rho&\wt\Phi_c\rix
=\mat{cc}\Omega&0\\0&\Omega_c\rix+z\mat{cc}\Delta_{11}&\Delta_{12}\\\Delta_{21}&\Delta_{22}\rix+O(z^2),
\]
where
\eq{factheta}
\mat{cc}\Delta_{11}&\Delta_{12}\\\Delta_{21}&\Delta_{22}\rix
=\mat{cc}M&0\\0&M_c\rix\mat{cc}\Omega \wt Q_1&\wt Q_1\\\Omega_c\wt Q_2&\wt Q_2\rix
\mat{cc}\wh B_{\rho-1,1}^{(\rho\rho)}&0\\0&\wh B_{\rho2}^{(\rho\rho)}\rix
\mat{cc}Q_1&Q_2\\Q_1\Omega&Q_2\Omega_c\rix.
\en
Following the standard perturbation theory, e.g., \cite{SteS90}, $\wh \Theta_\rho$ 
has an invariant subspace spanned by the columns of 
$$
\wh\Phi_\rho =\wt\Phi_\rho+z\wt \Phi_cY
$$
corresponding to the eigenvalues of $\Omega+z \Delta_{11}$
in the first order sense, where $Y$ satisfies
$$
\Omega_cY-Y\Omega+\Delta_{21}=0.
$$
Then
$$
(U+zE_U)F\wh \Phi_\rho (\Omega+z\Delta_{11}) 
= (V+zV_1)F\wh\Phi_\rho+O(z^2),
$$
and from which $\Omega_1=\Delta_{11}$.
We are able to express
$$
F\wh\Phi_\rho=\wt G_1 Q_1
+\Phi+O(z^2),
$$
where
\bstar
\wt  G_1&=&
\mat{ccccccc}0_{s_1\times s_\rho}^T&\ldots&
 0_{s_{\rho-1}\times s_\rho}^T&
\mat{c}I_{s_\rho}\\ 0_{(\rho-1)s_{\rho}\times s_\rho}\rix^T&
\mat{c}G_{\rho+1,\rho}\\0_{\rho s_{\rho+1}\times s_\rho}\rix^T&\ldots&
\mat{c}G_{k\rho}\\0_{(k-1)s_{k}\times s_\rho}\rix^T\rix^T,\\
\Phi&=& \mat{ccccccc}0_{s_1\times s_\rho}^T&\ldots& 
\Phi_{\rho-1}^T&\Phi_{\rho}^T&\Phi_{\rho+1}^T&\ldots&\Phi_{k}^T\rix^T,
\estar
and
\bstar
\Phi_{\rho-1}&=&z\mat{c}\wh B_{\rho-1,1}^{(\rho-1,\rho)}S_\rho^{-1}Q_1\Omega\\\vdots\\\wh B_{\rho-1,1}^{(\rho-1,\rho)}S_\rho^{-1}Q_1\Omega^{\rho-1}\rix,\quad
\Phi_{\rho}=\mat{c}zQ_2Y\\Q_1\Omega+zQ_2\Omega_cY\\\vdots
\\Q_1\Omega^{\rho-1}+zQ_2\Omega_c^{\rho-1}Y\rix,
\\
\Phi_{\rho+1}&=&
z\mat{c}C_{\rho+1}Q_1\Omega+G_{\rho+1,\rho}Q_2Y\\G_{\rho+1,\rho}Q_1\Omega
\\\vdots\\G_{\rho+1,\rho}Q_1\Omega^{\rho-1}\\(G_{\rho+1,\rho}S_\rho-\wh B_{\rho1}^{(\rho+1,\rho)})Q_1\rix,\quad
\Phi_{i}=
z\mat{c}C_iQ_1\Omega+G_{i\rho}Q_2Y\\G_{i\rho}Q_1\Omega\\0\\\vdots\\0\\-\wh B_{i-1,1}^{(i\rho)}Q_1\rix,\quad \rho+2\le i\le k.
\estar
Finally, let
$$
H:=(\Xi+\Xi_cP(t))RF\wh \Phi_{\rho}=\Xi RF\wh \Phi_{\rho}+O(t)=
H_0+t^{1/\rho}H_1+O(t^{2/\rho}).
$$
Then $H_0=\wt\Xi_\rho\mat{c}I_{s_\rho}\\G_\rho^{(\rho)}\rix Q_1$, which is just the 
constant term of $H$ in (\ref{ordh}), and
\bstar
H_1&=&\Xi_{\rho-1}\mat{c}\wh B_{\rho-1,\rho}^{(\rho-1,\rho)}S_\rho^{-1}Q_1\Omega
\\0_{(\rho-2)s_{\rho-1}\times r}\rix+\Xi_\rho
\mat{c}Q_2Y\\Q_1\Omega\\0_{(\rho-2)s_\rho\times r}\rix+\sum_{i=\rho+1}^k
\Xi_i\mat{c}C_iQ_1\Omega+G_{i\rho}Q_2Y\\G_{i\rho}Q_1\Omega\\0_{(k-2)s_k\times r}\rix\\
&=&\wt\Xi_{\rho-1}\mat{c}\wh B_{\rho-1,\rho}^{(\rho-1,\rho)}S_\rho^{-1}Q_1\Omega\\
Q_2Y\\C_{\rho+1}Q_1\Omega+G_{\rho+1,\rho}Q_2Y\\\vdots\\
C_kQ_1\Omega+G_{k\rho}Q_2Y\rix
+\wt\Xi_{\rho2}\mat{c}Q_1\Omega\\
G_{\rho+1,\rho}Q_1\Omega\\\vdots\\G_{k\rho}Q_1\Omega\rix\\
&=&\wt \Xi_{\rho-1}\mat{c}\wh B_{\rho-1,\rho}^{(\rho-1,\rho)}S_\rho^{-1}\\
0\\C_{\rho+1}\\\vdots\\C_k\rix Q_1\Omega+\wt\Xi_{\rho}\mat{c}I_{s_\rho}\\G_\rho^{(\rho)}\rix Q_2Y
+\wt\Xi_{\rho2}
\mat{c}I_{s_\rho}\\G_\rho^{(\rho)}\rix Q_1\Omega,
\estar
where 
\[
\wt\Xi_{\rho2}=\mat{cccc}\Xi_{\rho2}&\Xi_{\rho+1,2}&\ldots&\Xi_{k2}\rix,
\]
where $\Xi_{i2}$, $i=\rho,\ldots,k$, are given in (\ref{bkxi}).

We have proved the following theorem.
\begin{theorem}\label{thm1st}
Under the conditions of Theorem~\ref{thm2}, one has 
$$
(A+tD)(H_0+t^{1/\rho}H_1)=(H_0+t^{1/\rho}H_1)
(\lambda_0I+t^{1/\rho}\Omega+t^{2/\rho}\Delta_{11})+O(t^{2/\rho})
$$ 
where $\Delta_{11}$ is given in (\ref{factheta}), and $H_0, H_1$ are given above.
\end{theorem}
Note that by using $AH_0=\lambda_0H_0$, the above relation can be simplified to
$$
AH_1=\lambda_0H_1+H_0\Omega+O(t^{1/\rho}).
$$
Based on the formulas, the columns of $H_0$ are in the eigenvector subspace of $A$
spanned by the columns of $\wt \Xi_\rho$ and the columns of $H_1$ are in the
subspace spanned by the columns of $\wt\Xi_{\rho-1}$ and $\wt\Xi_{\rho2}$.
Note that the columns of $\wt\Xi_{\rho2}$ are the first generalized eigenvectors of $A$
corresponding to the Jordan blocks of sizes no less than $\rho$.

We should point out that although the eigenvalues of $\Omega$ are different from the
rest eigenvalues of $\Theta_\rho$, the matrix $H_0$ and even $H_0+t^{1/\rho}H_1$
 can be rank deficient, although  $H$ given in Theorem~\ref{thm2} has full column rank. 
 This is quite different from the standard first order perturbation results. 
\begin{example}\label{ex1}\rm
Consider the matrix
$$
A+tD=\mat{cccc}0&1&0&0\\0&0&1&0\\0&0&0&1\\t&0&0&0\rix,\qquad t>0.
$$
It has four eigenvalues
 $\lambda_j=t^{\frac14}e^{\frac12 i(j-1)\pi}$, $j=1,2,3,4$.
For 
$$
H=\mat{ccc}1&1&1\\\lambda_1&\lambda_2&\lambda_3\\
\lambda_1^2&\lambda_2^2&\lambda_3^2\\\lambda_1^3&\lambda_2^3&\lambda_3^3\rix
=\mat{cccc}1&&&\\&t^{1/4}&&\\&&t^{2/4}&\\&&&t^{3/4}\rix
\mat{crr}1&1&1\\1&i&-1\\1&-1&1\\1&-i&-1\rix,
$$
one has
$$
(A+tD)H=H\diag(\lambda_1,\lambda_2,\lambda_3).
$$
We have 
$$H_0=\mat{ccc}1&1&1\\0&0&0\\0&0&0\\0&0&0\rix,\qquad
H_1=\mat{ccc}0&0&0\\1&i&-1\\0&0&0\\0&0&0\rix.
$$
Although $H$ has full column rank
for any $t>0$, both $H_0$ and $H_0+t^{\frac14}H_1$
are rank deficient. 
\end{example}

\medskip
In the following, we consider the special case when $S_\rho$ has a semi-simple
eigenvalue $\gamma^{(\rho)}$ with algebraic multiplicity $r$. Suppose $\mu^{(\rho)}$
is a $\rho$th root of $\gamma^{(\rho)}$ and we choose $\Omega =\mu^{(\rho)} I_r$.
Let $Q_1,\wt Q_1^T\in\setC^{s_\rho\times r}$ have full column rank and satisfy
$$
S_\rho Q_1=\gamma^{(\rho)} Q_1,\quad \wt Q_1 S_\rho=\gamma^{(\rho)} \wt Q_1,\qquad
\wt Q_1 Q_1=I_r.
$$
Then for
\bstar
\wt \Phi_\rho&=&\mat{cccc} Q_1^T&\mu^{(\rho)} Q_1^T&\ldots&(\mu^{(\rho)} )^{\rho-1}Q_1^T\rix^T,\\
\wt \Psi_\rho&=&\frac1{\rho(\mu^{(\rho)} )^{\rho-1}}
\mat{cccc}(\mu^{(\rho)})^{\rho-1}\wt Q_1&\ldots&\mu^{(\rho)} \wt Q_1&\wt Q_1\rix,
\estar
one has
$$
\Theta_\rho \wt \Phi_\rho=\mu^{(\rho)} \wt\Phi_\rho,\quad
\wt\Psi_\rho\Theta_\rho=\mu^{(\rho)} \wt\Psi_\rho,\qquad \wt \Psi_\rho\wt\Phi_\rho=I_r.
$$
In this case,
$$
\Delta_{11}=
\frac1{\rho (\mu^{(\rho)} )^{\rho-2}} \wt Q_1(\wh B_{\rho-1,1}^{(\rho\rho)}+\wh B_{\rho2}^{(\rho\rho)})Q_1,
\qquad
\Delta_{21}= M_c (\Omega_c\wt Q_2\wh B_{\rho-1,1}^{(\rho\rho)}Q_1
+\mu^{(\rho)} \wt Q_2\wh B_{\rho2}^{(\rho\rho)}Q_1),
$$
where $M_c,\wt Q_2$ are the same as that given in (\ref{moq}),
$\Omega_c$ is given in (\ref{q2}), and 
$$
Y=(\mu^{(\rho)} I-\Omega_c)^{-1}\Delta_{21}.
$$
We have the same $H_0$ as before and
$$
H_1=\mu^{(\rho)}\wt \Xi_{\rho-1}\mat{c}\wh B_{\rho-1,\rho}^{(\rho-1,\rho)}S_\rho^{-1}\\
0\\C_{\rho+1}\\\vdots\\C_k\rix Q_1+\wt\Xi_{\rho}\mat{c}I_{s_\rho}\\G_\rho^{(\rho)}\rix Q_2(\mu^{(\rho)}I-\Omega_c)^{-1}\Delta_{21}
+\mu^{(\rho)}\wt\Xi_{\rho2}
\mat{c}I_{s_\rho}\\G_\rho^{(\rho)}\rix Q_1,
$$
where $Q_2$ is the same as the one in (\ref{q2}).
The eigenvalues corresponding to the perturbed invariant subspace are those of
$$
\lambda_0I_r+t^{1/\rho}\mu^{(\rho)}  I_r+\frac{t^{2/\rho}}{\rho (\mu^{(\rho)} )^{\rho-2}} \wt Q_1(\wh B_{\rho-1,1}^{(\rho\rho)}+\wh B_{\rho2}^{(\rho\rho)})Q_1+O(t^{3/\rho}).
$$
When $r=1$, it reduces to the first fractional order  perturbation results for eigenvectors.

\medskip
In principle, one can obtain the coefficient matrices for an arbitrary higher fractional order
in the same way by incorporating  all the blocks of $E_V(z)$ up to the same fractional order 
in the analysis. 
%
\section{Conclusions}\label{con}
We established the perturbation results for  invariant subspaces corresponding to
the eigenvalues that are perturbed from a single eigenvalue of the original matrix
with Jordan blocks of different sizes. For a matrix of an invariant subspace, 
we provide the fractional orders for all its 
blocks in the generic 
case. The results generalize the existing eigenvalue and eigenvector perturbation
theory in \cite{Lid65,MorBO97} and extend  the standard invariant subspace
perturbation theory as well. The situation is complicated in the nongeneric case and
further work needs to be done.

\medskip
\noindent
{\bf Acknowledgement}
The author thanks the anonymous referee and the handling editor for their comments and suggestions that helped to improve the paper.


\end{document}